\documentclass{article}

\usepackage{graphicx}
\usepackage{amsmath}
\usepackage[version=4]{mhchem}
\usepackage{siunitx}
\usepackage{longtable,tabularx}
\usepackage{geometry}
\usepackage{amsmath} 
\usepackage{amsfonts}
\usepackage{amssymb} 
\usepackage{graphicx} 
\usepackage{bm}  
\usepackage{setspace}
\usepackage{color}  
\usepackage{epstopdf}
\usepackage{enumitem}
\usepackage{algorithm}  
\usepackage{algorithmic} 
\usepackage{subfig}
\usepackage{multirow}

\title{Real-Time Optimal Control for Irregular Asteroid Landings Using Deep Neural Networks}

\author{Lin Cheng \thanks{Postdoctor, School of Aerospace Engineering, Tsinghua University, 10084 Beijing, China; chenglin5580@tsinghua.edu.cn. Member AIAA.},  
Zhenbo Wang\thanks{Assistant Professor, Department of Mechanical, Aerospace, and Biomedical Engineering, University of Tennessee, Knoxville, TN 37996, USA; zwang124@utk.edu. Member AIAA.},
Yu Song\thanks{Ph.D. Candidate, School of Aerospace Engineering, Tsinghua University, 10084 Beijing, China; yumail2011@163.com. Member AIAA },
Fanghua Jiang\thanks{\textbf{Corresponing Author}, Associate Professor, School of Aerospace Engineering, Tsinghua University, 10084 Beijing, China; jiangfh@mail.tsinghua.edu.cn. Senior Member AIAA.}
}

\date{}

\begin{document}

\maketitle

\begin{abstract}

Precise soft landings on asteroids are central to many deep space missions for surface exploration and resource exploitation. To improve the autonomy and intelligence of landing control, a real-time optimal control approach is proposed using deep neural networks (DNN) for asteroid landing problems wherein the developed DNN-based landing controller is capable of steering the lander to a preselected landing site with high robustness to initial conditions. First, to significantly reduce the time consumption of gravity calculation, DNNs are used to approximate the irregular gravitational filed of the asteroid based on the samples from a polyhedral method. Then, an approximate indirect method is presented to solve the time-optimal landing problems with high computational efficiency by taking advantage of the designed gravity approximation method and a homotopy technique. Furthermore, five DNNs are developed to learn the functional relationship between the state and optimal actions obtained by the approximate indirect method, and the resulting DNNs can generate the optimal control instructions in real time because there is no longer need to solve the optimal landing problems onboard. Finally, a DNN-based landing controller composed of these five DNNs is devised to achieve the real-time optimal control for asteroid landings. Simulation results of the time-optimal landing for Eros are given to substantiate the effectiveness of these techniques and illustrate the real-time performance, control optimality, and robustness of the developed DNN-based optimal landing controller.     

\end{abstract}

\section{Introduction}
Small bodies in the solar system such as asteroids and comets contain high scientific values, and various exploration missions have been carried out for exploring these small bodies \cite{yang2015fuel}. Landing a spacecraft on an asteroid for collecting high-quality resolution data and soil samples is a crucial step in the asteroid sample return missions, such as the finished NEAR \cite{bell2002asteroid} and Hayabusa 1 \cite{norton2008field} missions, and the ongoing Hayabusa 2 \cite{lauretta2012overview} and OSIRIS-REx \cite{lauretta2012overview} missions. However, there is a big challenge for asteroid landing due to the remote distance, irregular shapes, rapidly variable gravitational fields, and low thrust of the vehicle \cite{lantoine20061}. Considering the effects of these disturbances and uncertainties, the landing guidance and control approaches with high autonomy and reliability should be studied. 

The majority of the existing powered landing systems use two separate and independently systems for guidance and control, wherein a guidance system generates optimal trajectories that connect the lander's current state with its target state, and a control system is in charge of tracking these trajectories by determining the thrust magnitudes and directions \cite{furfaro2018recurrent}. Referring to this nominal trajectory tracking pattern, there are different derivative versions which are mainly focused on the improvement of trajectory optimization and trajectory tracking techniques. In generally, the asteroid landing problem has been formulated as an optimal control problem (OCP) and traditionally solved by the direct and indirect methods. Benefiting from the development of computers in recent years, pseudo-spectral method \cite{hu2016desensitized} and convex optimization techniques \cite{pinson2015rapid,yang2017rapid} are two effective direct methods that transform the original trajectory optimization problems into nonlinear programming problems (NLPs). Based on the study in \cite{pinson2015rapid}, the convex optimization method was improved in \cite{yang2017rapid} for solving the time-optimal trajectories of asteroids with irregular gravity. Meanwhile, indirect methods are also effective approaches for trajectory optimization with guaranteed solution accuracy and optimality \cite{jiang2012practical}. Some practical techniques are employed to further improve the solution performance of the indirect methods for asteroid landing trajectory optimization, such as the homotopy approach \cite{yang2018fasthomotopy} and the costate normalization \cite{sanchez2018real}. Additionally, many studies have been focused on developing more precise and robust trajectory tracking techniques with closed-loop control methods, such as linear sliding-mode control \cite{xiangyu2004autonomous}, terminal sliding-mode control \cite{lan2014finite},  convex optimization-based control \cite{cui2017intelligent}, and model predictive control\cite{alandihallaj2017soft}.

However, multiple issues need to be resolved for the trajectory tracking approach in asteroid landing missions \cite{furfaro2018recurrent}. First, both the direct and indirect methods suffer some drawbacks, including unpredictable iterative process and computational time, no guarantee of convergence to the optimal or even feasible solutions, and requirement of good initial guesses \cite{wang2017constrained}. Second, high-fidelity gravity calculations are computationally intensive which extends the time-consumption of trajectory optimization by hundreds of times \cite{lantoine20061}. For this reason, researchers usually have to design a reference trajectory in advance at the expense of losing flight autonomy. Third, the separate design of the guidance and control systems can easily lead to actuator saturation \cite{furfaro2018recurrent}. Moreover, the control system also affects the finial performance of the guidance system which increases the complexity of system optimization. Accordingly, it is worth exploiting more effective landing approaches with higher autonomy and reliability and guaranteed flight optimality and real-time performance.   


In addition to improving the performance of the indirect and direct methods, new control paradigms were also explored to enhance the autonomy and intelligence of the control system. In recent years, theoretical advancement in machine learning and improvement in artificial intelligence provide the possibility of real-time optimal control \cite{mnih2016asynchronous,gaudet2018deep}. According to the learning mode, machine learning algorithms are divided into three categories: supervised learning, unsupervised learning, and reinforcement learning \cite{goodfellow2016deep}. An improved reinforcement learning algorithm was proposed for the trajectory optimization of interplanetary solar sail spacecraft where an evolutionary algorithm was used to optimize the neural networks \cite{dachwald2004optimization}. A new actor-critic-identifier architecture was  presented for approximate optimal control of certain nonlinear systems \cite{bhasin2013novel}. Recently, an integrated guidance and control algorithm was developed for planetary powered descent and landing by applying the principles of reinforcement learning theory, and the trained policy could directly predict the thrust commands based on the lander's estimated state \cite{gaudet2018deep}. However, the current reinforcement learning algorithms are faced with two shortcomings: poor convergence property and hard reward function design. To address these issues, some researchers have begun to focus on combining the advantages of traditional trajectory optimization methods with the deep neural network (DNN) technology. Neural networks were used to approximate the value functions of the Hamilton$-$Jacobi$-$Bellman equations in \cite{medagam2009optimal,todorov2004optimality,lewis2003hamilton}. In addition, neural networks were also trained with supervised signals to approximate the optimal control of some simple motion tasks in \cite{effati2013optimal,effati2013optimal}. For example, the DNNs were trained in a supervised manner in \cite{sanchez2018real} using the optimal state-action pairs obtained via an indirect method, and the trained DNNs could drive the on-board landing systems with optimal control. Recurrent neural networks were also explored in a similar context in \cite{furfaro2018recurrent}. Recently, an interactive network learning algorithm was proposed to achieve real-time optimal control for solar sail orbit transfer missions with high terminal guidance accuracy \cite{chenglin2018TAES}. 

Different from the above mentioned works, a novel real-time optimal control approach to asteroid landings is investigated in this study by developing a DNN-based gravity model and proposing a new DNN learning architecture for landing control. The contributions of this paper are threefold: First, DNNs are developed to approximate the irregular gravitational field of the asteroid, and the well-trained DNNs can calculate the high-fidelity gravity with a much higher speed comparing to the polyhedral method \cite{werner1996exterior}. Second, incorporating the DNN-based gravity model in dynamics, an approximate indirect method is utilized to obtain optimal trajectories for landing, and a homotopy method is employed to improve the computational efficiency of the indirect method by supplying good initial guesses \cite{jiang2012practical}. With the help of the gravity approximation technique and the homotopy method, this approximate indirect method can quickly generate the optimal trajectories for asteroid landing missions with irregular gravitational fields. Third, new DNNs are developed to approximate the optimal control actions obtained by the approximate indirect method, and the resulting DNN-based landing controller can autonomously drive the lander from a random initial condition to the target point in a real-time optimal manner.         

This paper is organized as follows: In Sec.~\ref{sec: Problem Formulation and  Solution Architecture}, the time-optimal asteroid landing problem is formulated as a two-point boundary value problem, and an indirect-method-based network learning architecture is presented to achieve real-time optimal control for this asteroid landing mission. In Sec.~\ref{sec: Gravity Approximation}, a DNN is designed to approximate the gravitational field of the asteroid. Using this DNN-based gravity model in the dynamics, an approximate indirect method with a homotopy technique is presented based on the Pontryagin's Minimum Principle in Sec.~\ref{sec: DNN-Based Real-Time Optimal for Landing}, and then a DNN-based landing controller is developed through learning the optimal actions obtained by the approximate indirect method. In Sec.~\ref{sec: Simulation and Results}, numerical simulations are given to evaluate the effectiveness of the proposed techniques and illustrate the performance of the designed DNN-based landing controller. The main work of this study is summarized in Sec.~\ref{sec: conclusion}.

\section{Problem Formulation and Solution Architecture}
\label{sec: Problem Formulation and  Solution Architecture}

\subsection{Problem Formulation}

In this paper, we consider a time-optimal trajectory design problem for powered descent and landing on the surface of an asteroid. To achieve a soft pinpoint landing, a 3-dimensional (3D) \textbf{time-optimal trajectory control problem} is considered to find the thrust program that minimizes the following cost function
\begin{equation}
\label{equ: time-optimal  index}
 \mathop{\min}_{\bm{\mu, \alpha}}  \, t_f = \int_{t_0}^{t_f} 1 dt
\end{equation}
subject to the equations of motion \cite{yang2017rapid}
\begin{equation}
\label{equ: EoM }
\begin{array}{l}
 \dot{\bm{r}} = \bm{v}  \\
 \dot{\bm{v}} + 2 \bm{\Omega} \times \bm{v}+ \bm{\Omega} \times (\bm{\Omega} \times \bm{r}) = \bm{G}(\bm{r}) + \frac{T_{\max}\mu\bm{\alpha}}{m}  \\
 \dot { m } = - \frac { T_{\max} \mu } { I_{\mathrm {sp} }  g_{0} }
 \end{array}
\end{equation}
and the boundary conditions
\begin{gather}
 \bm { r } ( 0 ) = \boldsymbol { r } _ { 0 } , \quad \boldsymbol { v } ( 0 ) = \boldsymbol { v } _ { 0 } , \quad m ( 0 ) = m _ { 0 } \label{equ: initial state} \\
 \boldsymbol { r } \left( t _ { f } \right) = \boldsymbol {r}_{f} , \quad \boldsymbol { v } \left( t _ { f } \right) = \mathbf { 0 } \label{equ: terminal state} 
\end{gather}
where $\bm{r}=[x,y,z]^T$ and $\bm{v}=[v_x,v_y,v_z]^T$ are the lander's position and velocity respectively with respect to an asteroid-fixed coordinate frame, and $\bm{\Omega} = [0,0, \omega]^T$ denotes the angular velocity vector of the asteroid. The notation $\bm{G}$ represents the force vector due to the asteroid's gravitational field, $T_{max}$ is the maximum magnitude of the thrust, $\mu \in[0,1]$ is the thrust ratio, $\bm{\alpha}=[\alpha_1, \alpha_2, \alpha_3]$ is the direction vector of the thrust, and $m$ is the mass of the lander. The symbol $I_{\mathrm {sp}}$ denotes the thrust specific impulse, and $g_{0} =9.80665 {\rm m/s^2}$ is the standard acceleration of gravity at the Earth sea level. The landing site is assumed to be preselected, and the landing velocity is assumed to zero in the asteroid-fixed frame for the purpose of soft landing. Different from the traditional methods, the initial state variables $[\bm{r}_0,\bm{v}_0,m_0]^T$ of the lander are assumed to be unknown beforehand. That is to say, the landing controller developed in this study is expected to have autonomous landing capability.

The considered time-optimal problem of asteroid landing is an optimal control problem that was traditionally solved by the indirect and direct methods \cite{lantoine20061}. Indirect and direct methods show excellent performance on constraints satisfaction, optimality, and convergence property on common trajectory optimization tasks; however, it is extremely difficult for on-board generation of asteroid optimal landing trajectories \cite{furfaro2018recurrent}, due to the following reasons: First, solving the TPBVPs by indirect methods and the NLPs by direct methods needs good initial guesses \cite{wang2018minimum}. Second, there is no guarantee of convergence to an optimal or even feasible solution \cite{yang2018fasthomotopy}. Third, because extremely irregular and diverse asteroid shapes are involved, complicated gravitational field needs to be addressed. However, an accurate gravity calculation is computationally intensive, which leads to high time cost with hundreds of trajectory iterations \cite{lantoine20061}. In this study, DNNs accompanied by an improved indirect method are used to address these issues.

Considering the intensive computation time cased by the irregular gravitational filed, the time-optimal control problem is connected with a \textbf{gravity-free time-optimal control problem} with the same performance index and constraints, but the gravity $\bm{G}(r)$ is zero in the dynamical equations, that is
\begin{equation}
\label{equ: gravity-free EoM}
 \dot{\bm{v}} + 2 \bm{\Omega} \times \bm{v}+ \bm{\Omega} \times (\bm{\Omega} \times \bm{r}) = \frac{T_{\max}\mu\bm{\alpha}}{m} \\
\end{equation}
Since the time-consuming gravity calculation is no longer needed, the relatively easier problem for gravity-free time-optimal control can be quickly solved and its solutions could serve as good initial guesses for the original time-optimal control problem with the help of a homotopy technique \cite{jiang2012practical,yang2018fasthomotopy}.

\subsection{Overall Network Learning Architecture}

\begin{figure}[htbp]
\center
\vspace {2mm}
\centerline{\includegraphics[width= 4.5 in]{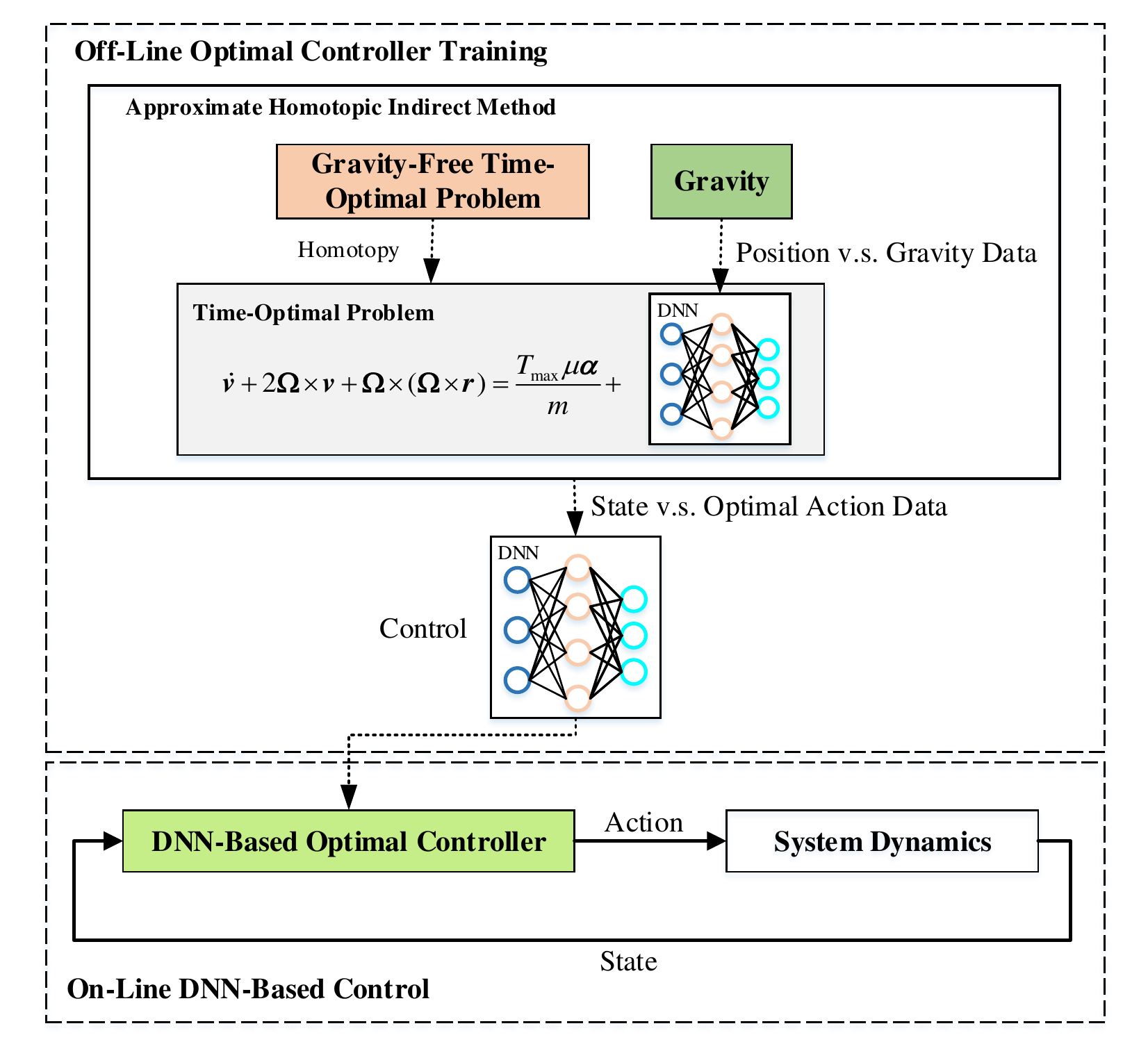}}
\vspace {-1mm}
\caption{The overall network learning architecture}
\label{fig:Actor-Indirect Method Networking Learning Architecture}
\vspace {-3mm}
\end{figure}

Taking advantage of high computational efficiency of the indirect methods and strong learning ability of the DNN technique, an indirect-method-based network learning architecture is proposed and depicted in Fig.~\ref{fig:Actor-Indirect Method Networking Learning Architecture} to achieve real-time optimal control for irregular asteroid landings. In this architecture, an indirect method is employed to solve the time-optimal landing problem and generate optimal state-action pairs. The control DNNs are trained off-line to approximate the functional relationship between the state and optimal actions, and can be used to generate optimal instructions for on-line landing control. Considering that the gravity calculation is usually computationally intensive \cite{lantoine20061}, a DNN is developed to learn the irregular gravitational field based on samples generated from a polyhedral method. The gravity approximation can significantly speed up the solution process of the indirect methods in generating optimal trajectories. Additionally, a homotopy approach in terms of the gravity is used to connect the time-optimal control problem with the gravity-free time-optimal control problem. Since the gravity-free time-optimal control problem does not contain the gravity $\bm{G}(\bm{r})$, it could be solved quickly. Thereafter, the original time-optimal control problem could be solved with the good initial costates provided by the gravity-free time-optimal control problem. With the gravity approximation and homotopy techniques, the resulting approximate indirect method can efficiently generate optimal landing trajectories which contain optimal actions corresponding to different flight states. Through learning the samples extracted from the obtained optimal trajectories, the DNN-based controller could establish a functional relationship between the states and optimal actions and enable optimal action predictions. Since it no longer needs to solve the time-consuming trajectory optimization problem onboard, the DNN-based landing controller enjoys the ability of real-time optimal control for asteroid landing missions.

\section{DNN Approximation of Gravity}
\label{sec: Gravity Approximation}

To address the issue that traditional gravity models are usually computationally intensive, a DNN is developed in this section to approximate the irregular gravitational field of the asteroid. This DNN-based gravity model is aimed to improve the calculation speed of gravity prediction with guaranteed accuracy.

\subsection{Gravitational Attraction Model}
The accurate irregular gravity is traditionally computed by two methods: the spherical harmonic expansion method and the polyhedral method \cite{werner1996exterior}. Although the spherical harmonic expansion method is easy to use, its drawback is that severe divergence appears close to the asteroid's surface. In contrast, considering the exact potential of the polyhedron anywhere of an asteroid in space, the polyhedron potential method is a high-fidelity approximation of the gravitational field in the region close the surface on an asteroid \cite{werner1996exterior}. Thus, the polyhedral method is employed to generate samples for DNN training in this study.

The exterior gravitational potential of a constant-density polyhedron can be derived analytically as follows \cite{werner1996exterior}

\begin{equation}
\label{equ: exterior gravitational potential}
\bm{G} = - G_u \rho\sum_{e\in edges} L_e \bm{E}_e \cdot \bm{r}_e + G_u \rho\sum_{f\in faces} \omega_f  \bm{F}_f \cdot \bm{r}_f
\end{equation}
where $G_u$ is the gravitational constant and $\rho$ is the density of the asteroid. $\bm{r}_e$ is a vector from the field point to an arbitrary point on each edge, $\bm{E}_e$ is a dyad defined in terms of the face and edge normal vectors associated with each edge, $L_e$ is a logarithmic term expressing the potential of a one-dimenstional straight wire, $\bm{r}_f$ is a vector from the field point to an arbitrary point on each face, $\bm{F}_f$ is the outer product of the face normal vectors, $\omega_f$ is the solid angle subtended by a face when viewed from the field point.

As we can see that the formula in Eq.~\eqref{equ: exterior gravitational potential} involves the computation of sums over all the edges and faces of the polyhedron. Since an accurate shape model requires a large number of faces and edges, this approach is computationally intensive, particularly when the indirect method runs tens of thousands of times to generate optimal trajectories. To address this issue, a DNN is developed to learn the gravitational force filed. Then, the trained network can replace the polyhedral item in the dynamical model, and the computational speed of predicting the gravity can be significantly improved. In the subsequent section, we will present the detailed design of the DNN. 

\subsection{DNN Approximation}

In this subsection, the detailed implementations of the gravity DNN are shown which includes the training data generation, data normalization, network architecture design, and optimization algorithm selection for network training. These algorithm implementations involve a series of scheme comparisons and parameter optimization, and due to the limited space, only the optimized results are provided here.

\begin{figure}[htbp]
\center
\vspace {2mm}
\centerline{\includegraphics[width= 4 in]{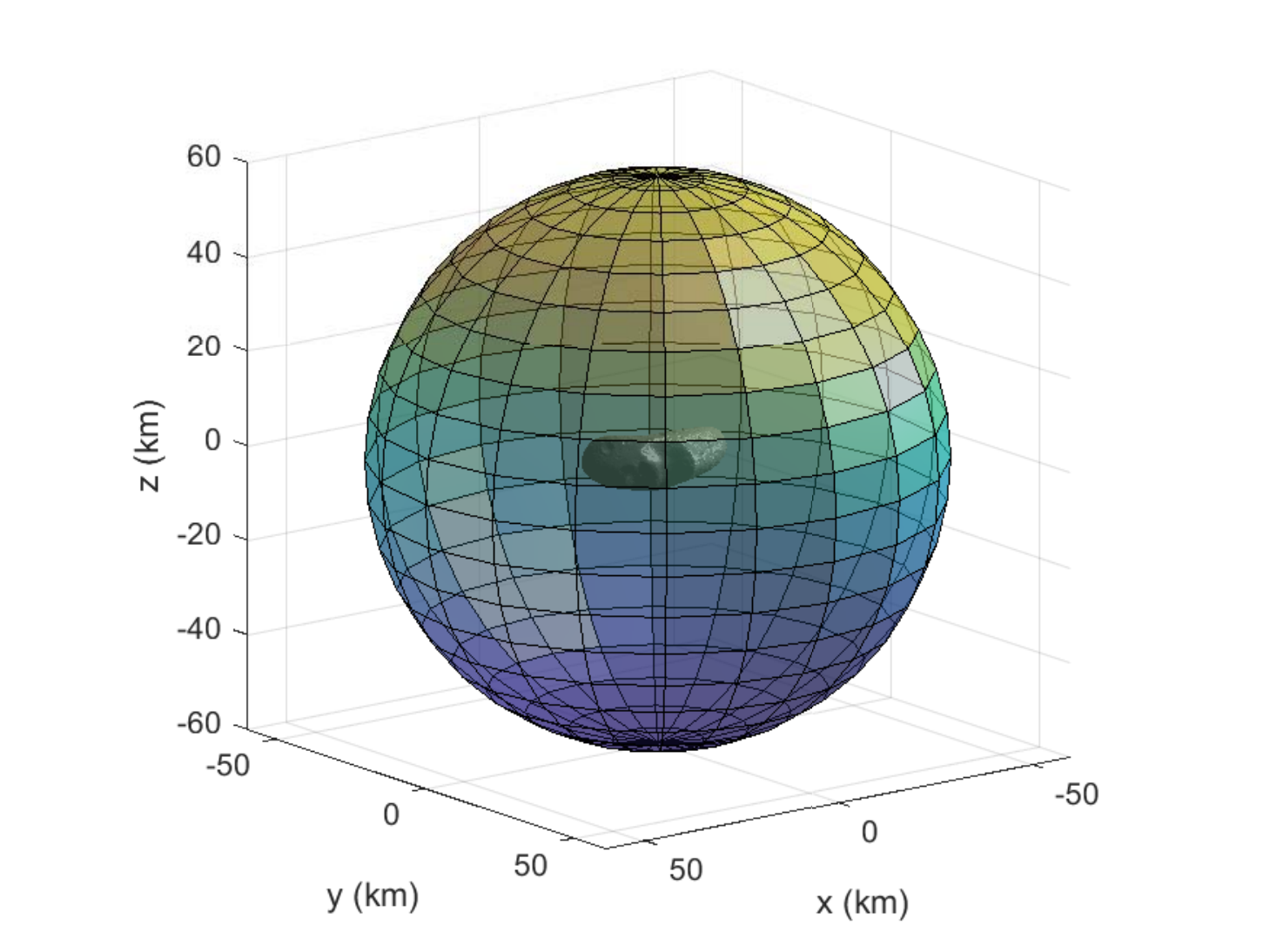}}
\vspace {-1mm}
\caption{Illustration of the gravity sampling}
\label{fig:Illustration of the Gravity Sampling}
\vspace {-3mm}
\end{figure}

The first step is to generate the training data of the gravitational field for DNN learning. This study is focused on the asteroid landing missions, hence, the developed gravity DNN mainly approximates the position-based gravity around the target asteroid. Taking the near-Earth asteroid 443 Eros as an example, it is a typical irregular-shaped asteroid with estimated diameters of $34.4\times11.2\times11.2$ km. Considering that the lander flies near and outside the target asteroid, the polyhedral method is used to calculate the gravity force vector in a spherical shell space where the radius is greater than 3 km and less than 60 km, as displayed in Fig.~\ref{fig:Illustration of the Gravity Sampling}. The position expressed in the spherical coordinate system is randomly chosen for the gravity calculation, and a total of 1,000,000 samples (position v.s. gravity force) are generated. According to the practical requirement of training, the samples are further randomly divided into two subsets with a ratio of 8:2 between the training dataset and test dataset. As their names imply, the training dataset is used to train the neural network, and the test dataset is used to evaluate the performance of network training and adjust the training parameters accordingly. Additionally, since the magnitude of the gravitational force is very small, we convert the input and output of the generated samples into the range of $[-1, 1]$. Simulation experiments have proved that this normalization technique can improve the data identification ability of the DNN and thus improve the learning effect of gravity approximation. In addition, there is a data anti-normalization process to obtain the actual gravity prediction from the trained DNN.

\begin{table}[htb]
\centering  
\caption{ The parameter setting of the gravity DNN for 443 Eros}
\label{tab:The architecture parameter of the gravity DNN}
\begin{tabular}{c|cccc}
 \hline
 \hline
                    &  Hidden Layer       & Output Layer  \\
 \hline
Activation Function &            ReLU               & Linear   \\
 \hline
Size                &         \multicolumn{2}{c}{ 6 Layers/256 Units }      \\
\hline
\hline
LR                &             \multicolumn{2}{c}{0.0001}   \\
\hline
N               &             \multicolumn{2}{c}{1000}   \\
\hline
\hline
\end{tabular}
\end{table}

A DNN is developed to learn the nonlinear functional relationship between the position and the gravitational force in this work. Since the gravity is dependent on the current position, deep feed-forward fully connected networks are employed. The design of the DNN architecture includes its activation function selection and size determination. The activation functions should be selected for both the output layers and the hidden layers. The activation functions of the output layers are determined by the value range of the sample output. Considering that the output of samples has no specific range limitation, thus, the gravity DNN adopts the Linear $[-\infty, \infty]$ activation function. Meanwhile, ReLU $[0,1]$ and Tanh $[-1,1]$ are two common activation functions for the hidden layers. In this approximation task, the ReLU activation function is found to be competent for the hidden layers and is employed in this study. The network size is another important user-defined parameter that affects the leaning effect, and both the depth and width should be carefully determined. Because the neuron numbers of the input and output layers are determined by the dimensions of the sample input and output, we only need to optimize the size of the hidden layers. As displayed in Table \ref{tab:The architecture parameter of the gravity DNN}, 6 layers and 256 units are found to be competent for gravity approximation of the asteroid 443 Eros.

The parameter vector $\bm{\omega}$ of the DNN comprises the weights and the biases of the neurons. By adjusting $\bm{\omega}$ in a way that the mean square error (MSE) between the DNN's predictions and the correct outputs from the training datasets is minimized, the DNN can accurately learn and approximate the gravitational field. The MSE is computed as follows
\begin{equation}
\label{equ: the mean square error of gravity DNN}
{L(\bm{\omega})} = \frac{1}{N}\sum\limits_{i = 1}^N {{{\left\|Net_{\mathrm {G}}(\bm{r}_i|\bm{\omega}) - {\bm{G}_i}\right\|^2}}} 
\end{equation}
where $Net_{\mathrm {G}}$ represents the gravity DNN and $\bm{G}_i$ is the correct gravity vector calculated by the polyhedral method for the ith sample. The symbol $N$ defines the number of samples that are randomly chosen from the training data to calculate the MSE in each training iteration. The algorithm Adam \cite{kingma2014adam} is employed to minimize the MSE. As displayed in Table \ref{tab:The architecture parameter of the gravity DNN}, the batch number $N$ is chosen as 1000 and the learning rate (LR) is set to be 0.0001.

Finally, the pseudo code of the training algorithm is summarized in Algorithm \ref{code: Network interactive training algorithm}, which is implemented in Python environment on a typical desktop using Tensorflow \cite{abadi2016tensorflow}, and the symbol $M$ denotes the maximum iteration number for the network training. To summarize, a specific DNN is developed in this subsection to approximate the gravitational field of our selected asteroid, and we will evaluate the approximation performance of the network $Net_{\mathrm {G}}$ through multiple experiments in the following subsection.

\begin{algorithm}[htb]  
\caption{The DNN-based gravity approximation}  
\label{code: Network interactive training algorithm}  
\begin{algorithmic}[1]  
\STATE Randomly initialize $Net_{\mathrm {G}}(\bm{r}|\bm{w})$ with weights $\bm{w}$	
\FOR{episode = 1, $M$}
\STATE   Select a random minibatch of $N$ samples $[\bm{r}, \bm{G}]$ from $\mathcal{R}$ \
\STATE   Update $Net_{\mathrm {G}}(\bm{r}|\bm{w})$ with Adam algorithm by minimizing the loss: \\
 $ {L(\bm{\omega})} = \frac{1}{N}\sum\limits_{i = 1}^N {{{\left\|Net_{\mathrm {G}}(\bm{r}_i|\bm{\omega}) - {\bm{G}_i}\right\|^2}}} $ \\
 \ENDFOR   
\end{algorithmic}  
\end{algorithm}

\subsection{Accuracy Evaluation of DNN-based Gravity Model}

\begin{figure*}[htbp]
\center
\vspace {2mm}
\centerline{\includegraphics[width= 6 in]{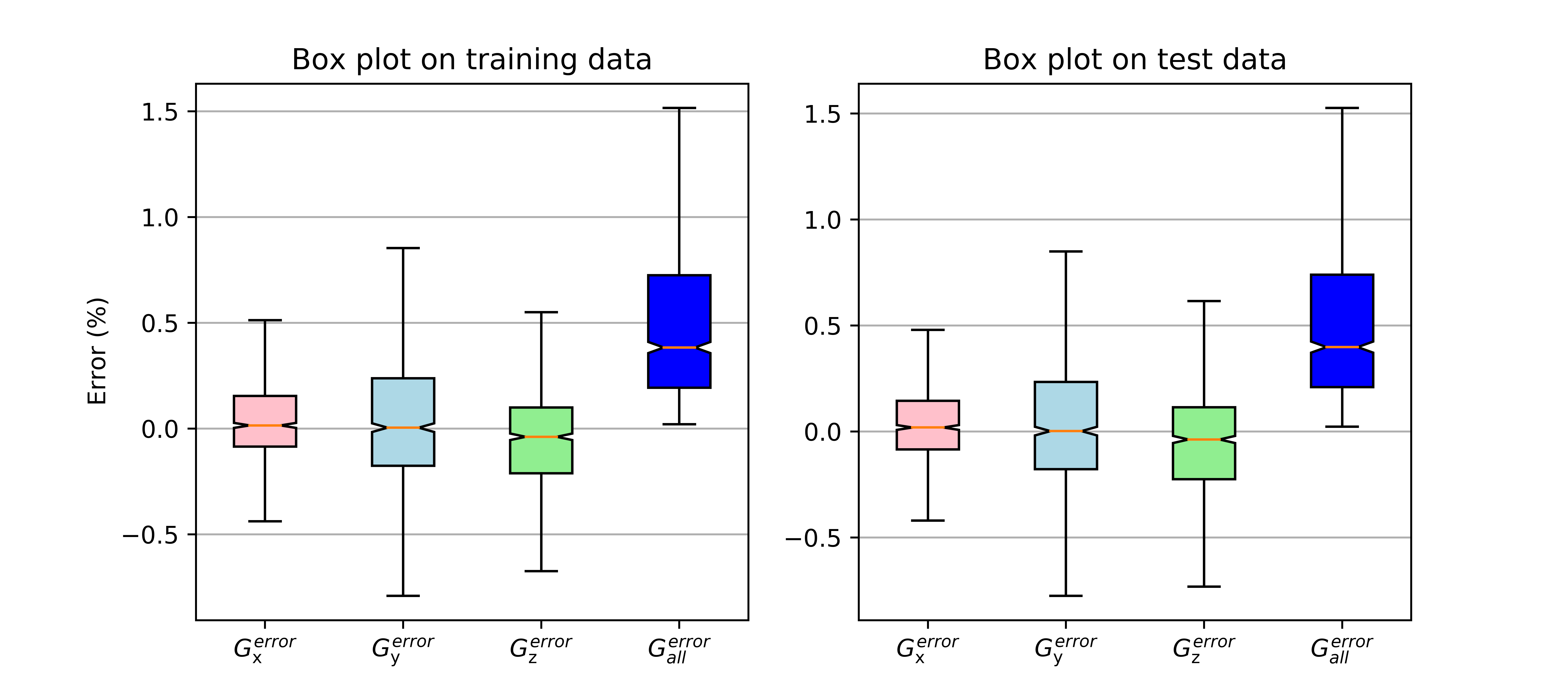}}
\vspace {-1mm}
\caption{Box plot of network $Net_{\mathrm {G}}(\bm{r}|\bm{w})$ prediction error}
\label{fig: Box plot of gravity network}
\vspace {-3mm}
\end{figure*}

Next, our attention turns to evaluate accuracy of the designed DNN for the irregular gravitational filed of the asteroid. The first experiment is to quantitatively analyze the approximation error of the gravity DNN where the relative errors projected onto the three axes and the synthetic relative error are defined as   
\begin{equation}
\label{equ: relative error vector}
\bm{G}_{\mathrm {axis}}^{\mathrm {error}} = \frac{Net_{\mathrm {G}}(\bm{r}|\bm{\omega})-{\bm{G}}}{\left\|{\bm{G}}\right\|_2} \times 100 \%
\end{equation}
and
\begin{equation}
\label{equ: sum of relative error }
\bm{G}_{\mathrm {all}}^{\mathrm {error}} = \frac{\left\|Net_{\mathrm {G}}(\bm{r}|\bm{\omega})-{\bm{G}}\right\|_2}{\left\|{\bm{G}}\right\|_2} \times 100 \%
\end{equation}
respectively, where the subscript ``${\mathrm {axis}}$'' represents one of the three axes x, y, and z, and $\Vert \cdot \Vert_2$ denotes the Euclidean norm of a vector. For the 443 Eros example, Fig.~\ref{fig: Box plot of gravity network} shows these four errors of 10,000 samples from the training dataset and test dataset using two box plots which display the distribution of four errors based on five number summary: minimum, first quartile, median, third quartile, and maximum. As can be seen from this figure, the relative errors along these three axes stay within the range of $\pm 0.75\%$, and most of their absolute values are less than $0.2\%$. Meanwhile, as a  combination of the three axial errors, the synthetic relative error is relatively bigger with a  median of about $0.3\%$ and a maximum of less than $1.6\%$. More than $3/4$ of the errors are less than $0.75\%$. Through this experiment, we can conclude that the developed gravity DNN can approximate the gravitational field of an asteroid with small prediction errors. However, considering that the landing is a long-lasting process, we should further evaluate the cumulative effect of the approximation error.

\begin{figure}[htbp]
\center
\vspace {2mm}
\centerline{\includegraphics[width= 6.2 in]{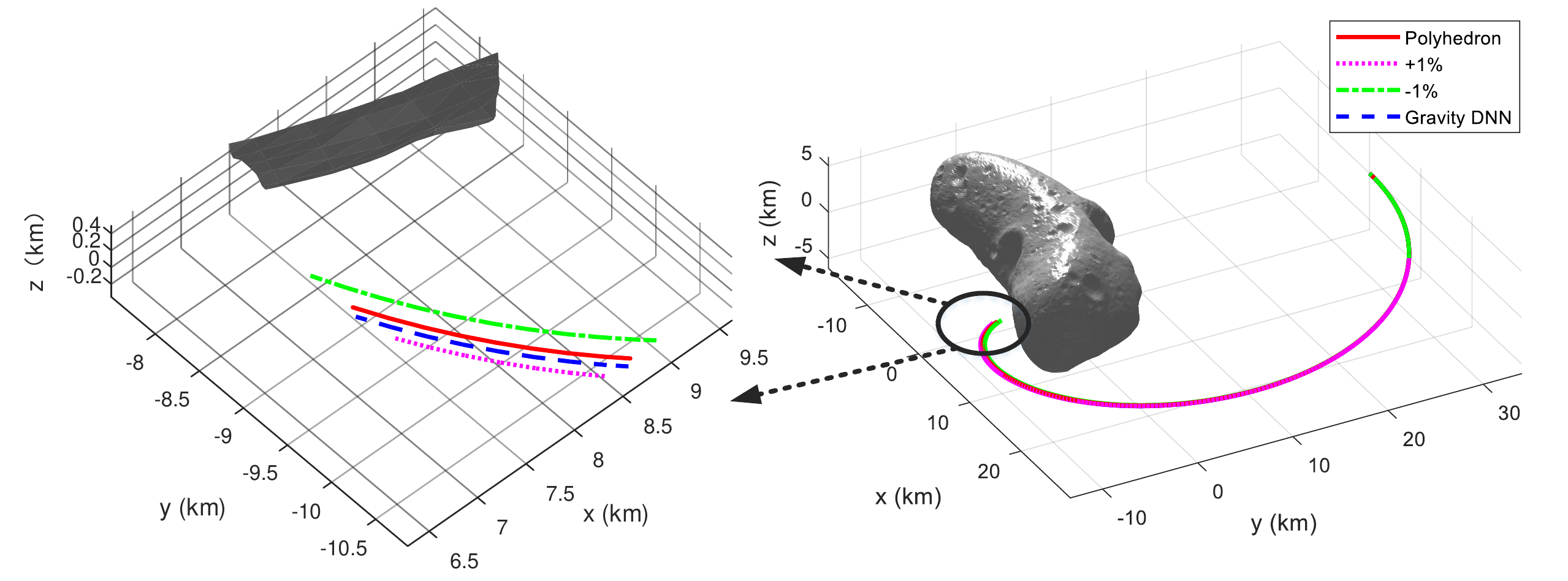}}
\vspace {-1mm}
\caption{Ballistic trajectory driven by different gravity models}
\label{fig: Ballistic trajectory driven by different gravity model}
\vspace {-3mm}
\end{figure}

Since the polyhedral method is also an approximate gravity model, and there must be approximation errors due to the effects of uneven distribution of mass, density error, polyhedron calculation error, etc. In our second experiment, we compare the orbit cumulative difference between the assumed deviation of the polyhedral model and the approximation error of the DNN. Taking the 443 Eros as an example, Fig.~\ref{fig: Ballistic trajectory driven by different gravity model} displays four trajectories of an unpowered lander that are driven by the ideal polyhedral model, the polyhedral model with $\pm1\%$ deviation, and the gravity DNN. The left subgraph is an enlarged view of the end of these four trajectories in the right subgraph. As it is shown in the figure that the trajectories forced by the gravity with errors would gradually deviate from the nominal trajectory. Since there are numerical difference between the assumed deviation of the polyhedral model and the approximation error of the DNN, the deviations of the trajectories are also different. We can see from the figure that the trajectory deviation due to the DNN approximation error is much smaller than the ones due to $\pm 1\%$ shifting of the polyhedral model. Considering that $\pm 1\%$ approximation error is usually acceptable in traditional gravitational models, we can conclude that the gravity DNN has outstanding approximation accuracy.


\begin{figure}[htbp]
\center
\begin{minipage}[t]{0.48\linewidth}
\centerline{\includegraphics[width= 3 in]{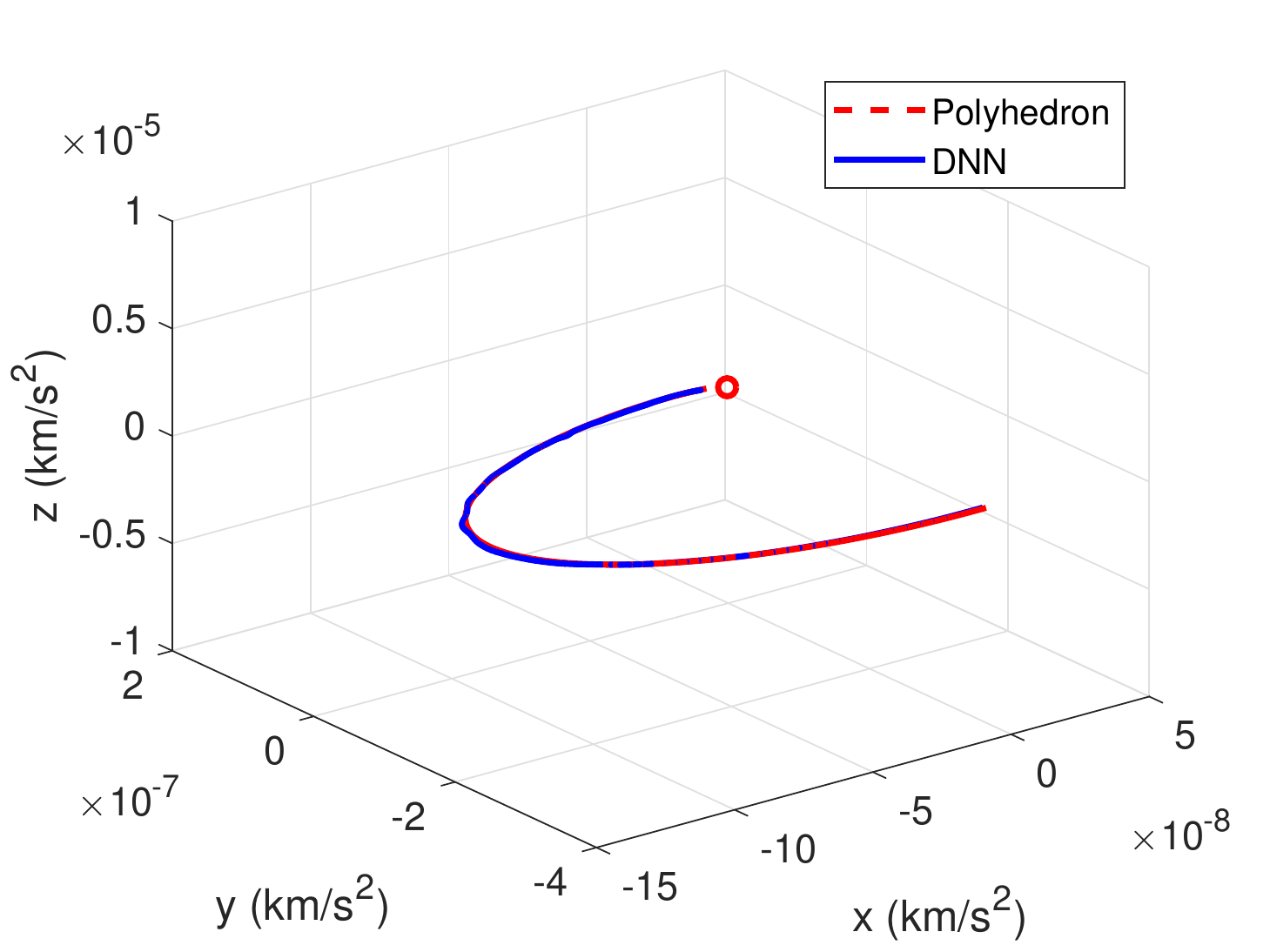}}
\centerline{\textbf{(a)}}
\end{minipage}
\begin{minipage}[t]{0.48\linewidth}
\centerline{\includegraphics[width= 3 in]{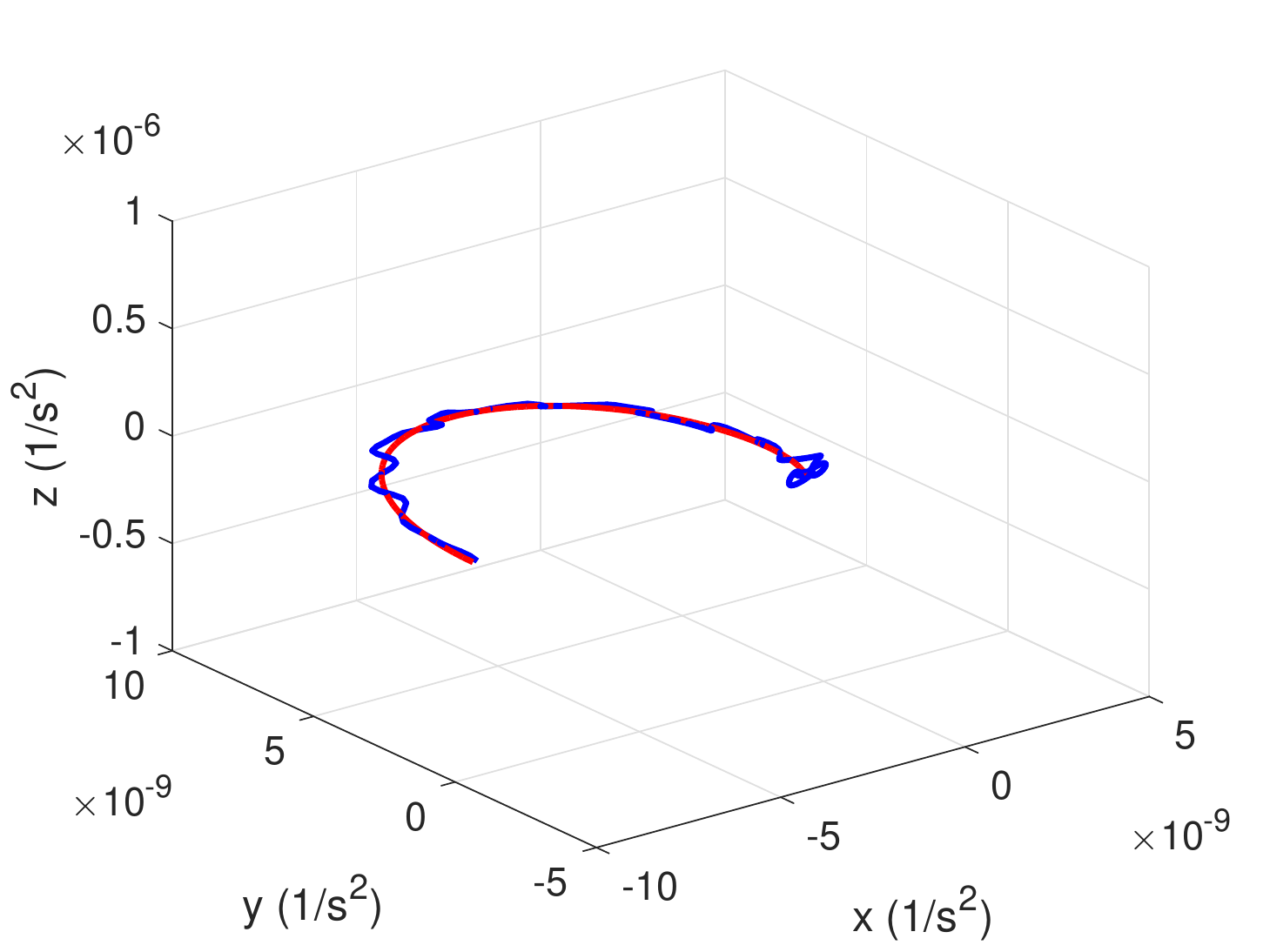}}
\centerline{\textbf{(b)}}
\end{minipage}
\begin{minipage}[t]{0.48\linewidth}
\centerline{\includegraphics[width= 3 in]{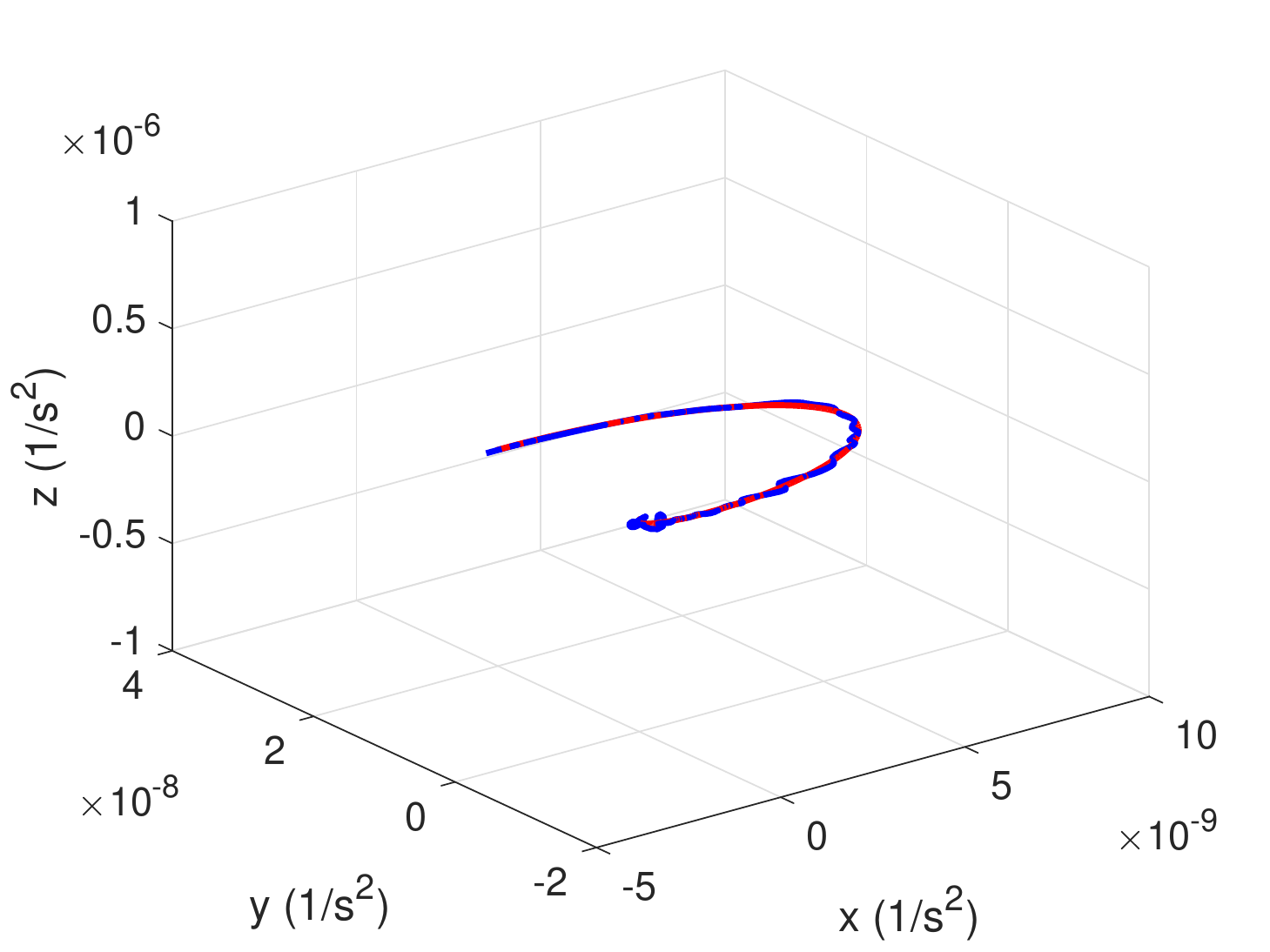}}
\centerline{\textbf{(c)}}
\end{minipage}
\begin{minipage}[t]{0.48\linewidth}
\centerline{\includegraphics[width= 3 in]{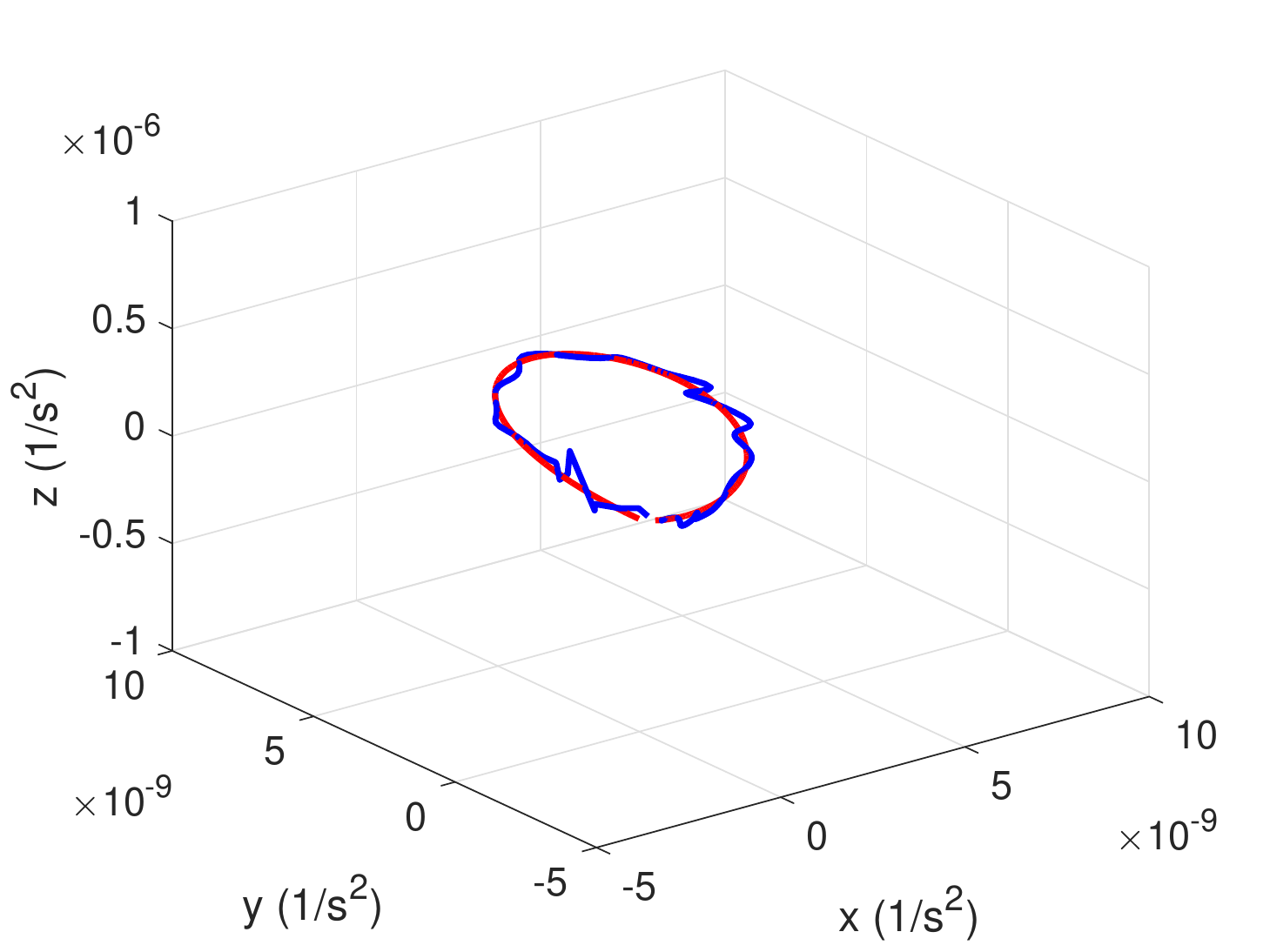}}
\centerline{\textbf{(d)}}
\end{minipage}
\caption { Approximation results of the DNN-based gravity model, where (a) displays its gravity predictions compared with the target outputs calculated by the polyhedral method, while (b), (c) and (d) illustrate the gravity predictions of the gravity partial derivatives}
\label{fig: GravityCompare}
\end{figure}

The third experiment is to show the approximation results of the gravity DNN using graphs in the  Cartesian coordinates. The gravitational force at one position of an asteroid is a 3D vector. Figure \ref{fig: GravityCompare} (a) displays the vector-end curve of the gravity in a randomly chosen unpowered flight trajectory. As can be seen from the figure, the predictions of the trained DNN are almost the same as the results of the polyhedral method. Additionally, another advantage of the DNN-based gravity model is that the DNN is an explicit function $Net_{\mathrm {G}}(\bm{r}|\bm{\omega})$, and the partial derivatives of the gravity prediction (DNN output) with respect to the position (DNN input) can be obtained explicitly. And also, the subsequent indirect method would use these partial derivatives. Considering that the polyhedral method needs extra calculation of the partial derivatives, the DNN-based gravity model has another advantage in reducing the computational time. Subgraphs \ref{fig: GravityCompare} (b)-(d) compare the vector-end curves of the partial derivatives $\partial \bm{G_x}/\partial \bm{r}$,$\partial \bm{G_y}/\partial \bm{r}$,$\partial \bm{G_z}/\partial \bm{r}$ under these two models. As can be seen from the subgraphs that the gravity DNN behaves good approximation ability to the partial derivatives of the gravity.

Through learning the samples of gravity versus position obtained by the polyhedral method, we obtain a new gravity model represented by DNNs. This DNN-based gravity model can quickly and accurately predict the gravity forces and the gradient of the gravitational field with respect to the position. The approximation accuracy of this DNN-based gravity model has been verified based on the above three experiments, and its real-time performance will be discussed in Section \ref{sec: Simulation and Results}. The design DNN-based gravity model plays an important role in generating optimal landing trajectories in the subsequent section, and the trajectory generation efficiency will be significantly improved due to the excellent real-time performance of this model.

\section{DNN-Based Real-Time Optimal Landing Control}
\label{sec: DNN-Based Real-Time Optimal for Landing}

Based on the DNN-based gravity model in Section \ref{sec: Gravity Approximation}, an approximate indirect method is developed in this section to generate optimal trajectories for asteroid landings with irregular gravitational fields. Then, DNNs are designed to approximate the optimal control actions with respect to the flight states through off-line learning samples extracted from the generated optimal trajectories. Since there is no need to solve optimal control problems onboard, the resulting DNN-based landing controller is capable of optimal control with great potential for real-time asteroid landing applications.  

\subsection{Approximate Indirect Method for Optimal Landing Trajectory}

In this subsection, an approximate indirect method is developed for the time-optimal landing problem based on the above DNN-based gravity model. Applying the Pontryagin's Minimum Principle, the Hamiltonian of the time-optimal problem is formulated as 

\begin{equation}
\label{equ: Hamiltonian}
H = \bm{\lambda}_r \cdot \bm{v} + \bm{\lambda_v} \cdot \left(- 2 \bm{\Omega} \times \bm{v} - \bm{\Omega} \times (\bm{\Omega} \times \bm{r})   + Net_{\mathrm {G}}(\bm{r}|\bm{\omega} \right) +  \frac{T_{max} \mu}{m} \bm{\alpha} ) + \lambda_m \cdot \left(- \frac{T_{max} \mu}{I_{sp}g_0}\right) + \lambda_0 
\end{equation}
where $\bm{G}(r)$ is replaced with the trained network $Net_{\mathrm {G}}(\bm{r}|\bm{\omega})$. By introducing a positive costate normalization factor $\lambda_0$, the initial costate vector $\bm{\lambda} = [\lambda_0, \bm{\lambda}_r, \bm{\lambda_v}, \lambda_m]$ can be restricted on the surface of a unit sphere and satisfies
\begin{equation}
\label{equ: constate vector nomralization equal 1}
\Vert \bm{\lambda}(t_0) \Vert_2 = 1
\end{equation}

To minimize the Hamiltonian $H$ in Eq.~\eqref{equ: Hamiltonian}, the optimal control law of thrust direction and magnitude can be obtained as 
\begin{equation}
\label{equ: optimal control thrust direction}
\begin{array}{l}
 \bm{\alpha}^* = -\frac{\bm{\lambda_v}}{\Vert\bm{\lambda_v}\Vert_2} \\
 \mu^* = 1
\end{array}  
\end{equation}
where $\Vert \cdot \Vert $ represents the Euclidean norm and the superscript $^*$ denotes the optimal control. It has been proved that the thrust of the lander is always maximized for the time-optimal landing missions \cite{Junfengli2007DynamicsandControl}.

The costate differential equations corresponding to the Euler-Lagrange equations are derived by $\dot{\bm{\lambda}} = - \partial{H} / \partial{\bm{x}}$ and shown below 
\begin{equation}
\label{equ: lambda_r dot}
\begin{array}{l}
\dot{\bm{\lambda}}_r  = - \frac{\partial{H}}{\partial{\bm{r}}} = \bm{\Omega} \times (\bm{\Omega} \times \bm{\lambda_v})  -  \frac{\partial{Net_{\mathrm {G}}(\bm{r}|\bm{\omega})}}{\partial{\bm{r}}}\cdot \bm{\lambda_v} \\
\dot{\bm{\lambda}}_v  = - \frac{\partial{H}}{\partial{\bm{v}}} = = -\bm{\lambda_r} - 2 \bm{\Omega} \times \bm{\lambda_v} \\
\dot{\bm{\lambda}}_m  =  - \frac{T_{max}\mu \Vert \bm{\lambda_v}\Vert }{m^2} 
\end{array}
\end{equation}
where ${\partial{Net_{\mathrm {G}}(\bm{r}|\bm{\omega})}}/{\partial{\bm{r}}}$ approximates the gradient vector of the gravitational force with respect to the corresponding position, and its value is calculated by a back-propagation algorithm in Tensorflow.

According to the static condition and the transversality condition, the optimal arrival time $t_f$ and the final mass costate are determined as
\begin{gather}
 H(t_f) = 0 \label{equ: terminal static condition} \\
 \lambda_m(t_f) = 0 \label{equ: terminal transversality condition} 
\end{gather}

Combining the above derivation, a two-point boundary value problem (TPBVP) is obtained wherein the corresponding shooting function is  
\begin{equation}
\label{equ: shooting equation}
 \Phi(\bm{z})= [\bm{r}(t_f)-\bm{r}_f, \bm{v}(t_f)-\bm{v}_f, \lambda_m(t_f), \Vert \bm{\lambda}(t_0) \Vert_2 - 1, H(t_f)] = 0 \\
\end{equation}
and the corresponding designed variable vector is $\bm{z}=[\lambda_0, \bm{\lambda}_r, \bm{\lambda}_v, \lambda_m, t_f]$. This TPBVP is a nonlinear multi-variable root-finding problem and MINPACK’s hybrd and hybrj routines (modified Powell method) \cite{more1980user} are found to be competent to solve this problem. 

Although the DNN-based gravity model consumes much smaller time compared with the polyhedral model, high-fidelity integration of Eq.~\eqref{equ: EoM } is still time consuming. Moreover, the convergence of the solution method for the TPBVP depends on a good initial costate vector, and the computational time increases when a large number of iterations are required. Since substantial optimal trajectories are needed to train the landing controller, a fast approach to providing good initial costates for the approximate indirect method is necessary. To this end, a homotopy approach in terms of the gravity is used to connect the time-optimal landing problem with the gravity-free time-optimal landing problem, which can be solved much more quickly because the time-consuming calculation of gravitational force is not needed as shown in Eq.\eqref{equ: gravity-free EoM}. Using a homotopy factor $\varepsilon$, we can rewrite the differential equation of the velocity as
\begin{equation}
\label{equ: EoM of gravity-free time-optimal landing problem}
\begin{array}{l}
 \dot{\bm{v}} + 2 \bm{\Omega} \times \bm{v}+ \bm{\Omega} \times (\bm{\Omega} \times \bm{r}) = \varepsilon \cdot Net_{\mathrm {G}}(\bm{r}|\bm{\omega}) + \frac{T_{\max}\mu\bm{\alpha}}{m}  
 \end{array}
\end{equation}
where the parameter  $\varepsilon$ links the gravity-free time-optimal landing problem ($\varepsilon=0$) with the time-optimal landing problem ($\varepsilon=1$). After obtaining the solution to the gravity-free time-optimal landing problem, we will gradually increase the parameter $\varepsilon$ following a certain series and take the solution obtained in the current iteration as an initial guess for the next iteration. Finally, the solution to the original time-optimal landing problem can be obtained. Simulations show that because the gravity forces of asteroids are usually very small, the solutions to the gravity-free time-optimal landing problem can always provide good initial costates for the time-optimal landing problems, i.e., the homotopy takes one step in most scenarios.

\subsection{Generation of Training Data}

\begin{figure}[htbp]
\center
\begin{minipage}[t]{0.48\linewidth}
\centerline{\includegraphics[width= 3 in]{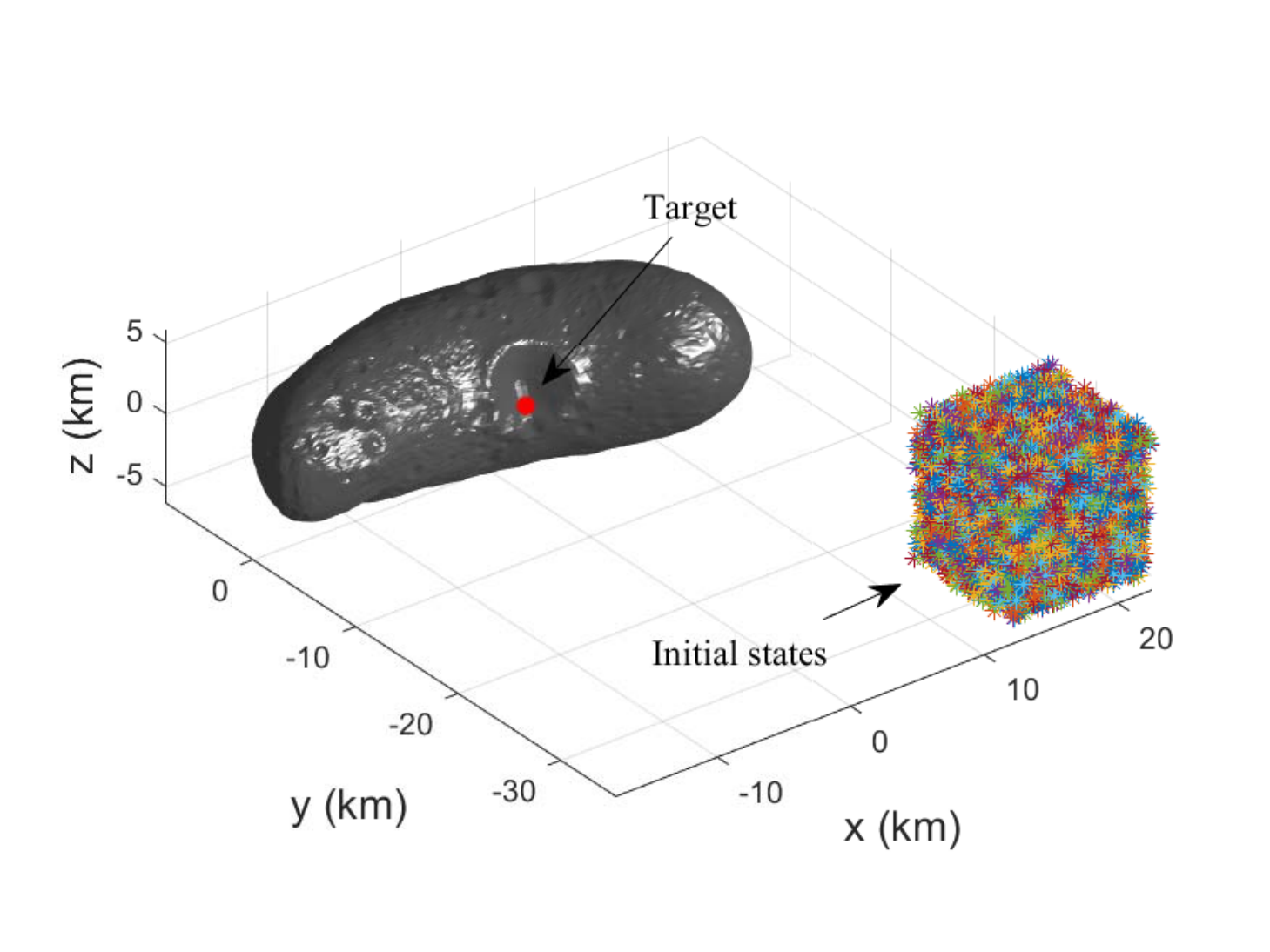}}
\centerline{\textbf{(a)}}
\end{minipage}
\begin{minipage}[t]{0.48\linewidth}
\centerline{\includegraphics[width= 3 in]{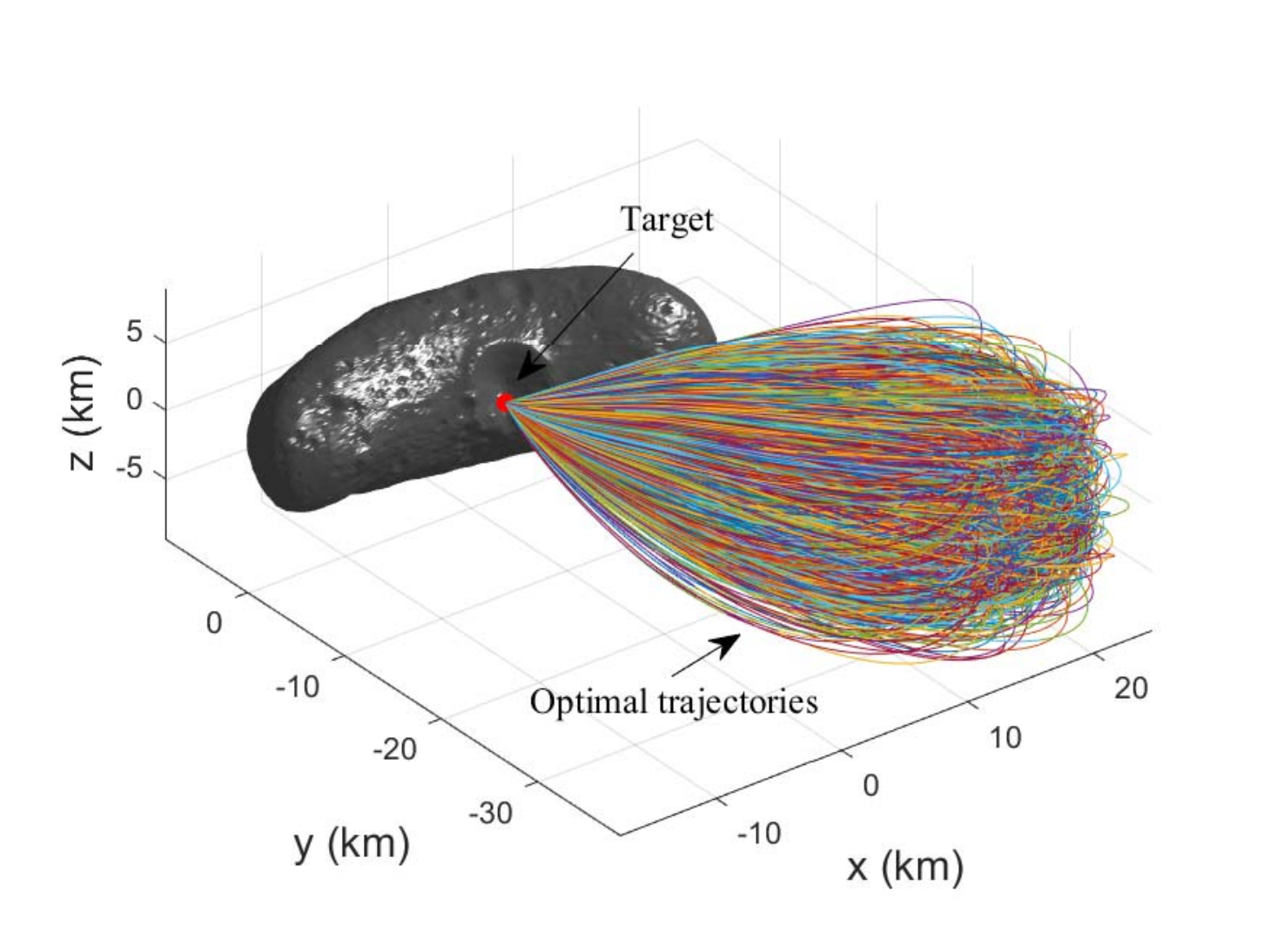}}
\centerline{\textbf{(b)}}
\end{minipage}
\caption{Data generation for the control DNN training where (a) shows the initial states of the lander and (b) displays the optimal trajectories obtained from the approximate indirect method }
\label{fig: DataGeneration}
\end{figure}

In traditional guidance and control approaches to asteroid landings, the lander's initial states are assumed to be known and the time-consuming trajectory planning can be conducted in advance. However, this assumption is almost impossible due to the  remoteness of the asteroids, rapidly variable gravitational fields, and uncertain force such as the solar wind and the aberration of light. In this study, we aim to develop a DNN-based intelligent landing controller which can steer the lander from an arbitrary initial state to the target point with high autonomy and reliability. To achieve this goal, we should get enough samples about the flight states versus optimal control actions and develop control DNNs to approximate the functional mapping relationship. In accordance with Bellman's optimality principle, knowing the optimized vector $\bm{z}$ satisfying $\Phi(\bm{z})$ is equivalent to knowing the whole optimal trajectory (referring to Eq.~\eqref{equ: shooting equation}). Consequently, any state-action solution pair $[\bm{x}, \bm{u}^*]$ along the optimal trajectory satisfies the first-order optimality condition, where $\bm{x}=[\bm{r}, \bm{v}]$ denotes the flight state and $\bm{u}^*$ consists of the corresponding optimal control $[\bm{\alpha}^*, \mu^*]$, time-to-go $t_f^*$, and expected mass consumption $\Delta m$. The trajectories generated by the approximate indirect method should cover the whole space that the landing process may happen. Taking the Eros as an example, we assume that the lander initiates its landing mission from a 3-dimensional cubic space as shown in Fig.~\ref{fig: DataGeneration} (a). The side length of this cube is 10 km, and the lander's initial position is randomly generated in this cube. Meanwhile, the lander's three initial velocity components are randomly chosen in the range of $[-0.01/1000, 0.01/1000] $ km/s. A total of 10,000 optimal trajectories (Fig.~\ref{fig: DataGeneration} (b)) randomly starting from this cubic space are obtained via the approximate indirect method, and for every trajectory, 200 optimal state-action pairs are extracted with equal time intervals. Then, the resulting 2,000,000 optimal state-action pairs are divided into two subsets for the DNNs: the training dataset with 80\% of the data and the test dataset with 20\% of the data.


\subsection{Network Optimization and Training}

\begin{table}[htb]
\centering  
\caption{ Parameters setting of the control DNNs}
\label{tab:The parameter setting of the control DNN}
\begin{tabular}{c|cccc}
 \hline
 \hline
                                  &                                               & Activation Function & Size                                 & Loss\\
 \hline
\multirow{2}{*}{$Net_{\alpha}$}   & Hidden Layer                                  & ReLU                & \multirow{2}{*}{6 Layer/512 Units  } & \multirow{2}{*}{ MAE } \\
\cline{2-3}
                                 
                                  & Output Layer                                  & Tanh                & \\
\hline
\multirow{2}{*}{$Net_{t_f}$}      & Hidden Layer                                  & ReLU                & \multirow{2}{*}{6 Layer/256 Units }  & \multirow{2}{*}{ MSE }\\
\cline{2-3}
                                  & Output Layer                                  & Softplus            & \\
\hline
\multirow{2}{*}{$Net_{\Delta m}$} & Hidden Layer                                  & ReLU                & \multirow{2}{*}{6 Layer/256 Units }  & \multirow{2}{*}{ MSE }\\
\cline{2-3}  
                                  & Output Layer                                  & Softplus            & \\
\hline
\hline
LR                                & \multicolumn{4}{c}{0.0001}                \\
\hline
N                                 & \multicolumn{4}{c}{1000}                   \\
\hline
\hline
\end{tabular}
\end{table}

Five DNNs, $Net_{\alpha1}$, $Net_{\alpha2}$, $Net_{\alpha3}$, $Net_{t_f}$, and $Net_{\Delta m}$, are designed to approximate the state-based optimal actions, which correspond to the three thrust directions ($\alpha_1^*$ , $\alpha_2^*$, and $\alpha_3^*$), time-to-go $t_f^*$, and mass consumption $\Delta m^*$. In order to avoid the approximation interference between control instructions, three independent DNNs, $Net_{\alpha1}$, $Net_{\alpha2}$, and $Net_{\alpha3}$ with the same structure are used to approximate the thrust directions. Meanwhile, there is no need to learn the thrust ratio $\mu^*$ since it is known to be a constant of one according to Eq.~\eqref{equ: optimal control thrust direction}. Table \ref{tab:The parameter setting of the control DNN} exhibits the structures and parameters of the five DNNs, including the selection of activation functions, network size determination, loss definition, and user-defined parameters. With regard to the selection of activation functions, all five DNNs adopt ReLU $[0, \infty]$ for the hidden layers. Meanwhile, the activation functions for the output layers are determined based on their output value ranges. Specifically, Tanh $[-1,1]$ is used for the three DNNs for the thrust directions, while Softplus $[0, \infty]$ is adopted in the DNNs of $Net_{t_f}$ and $Net_{\Delta m}$ since the time-to-go $t_f^*$ and mass consumption $\Delta m^*$ are always greater than zero. Loss functions are used to quantify the difference between the predictions of the DNN and the target outputs from the training data. In this regard, the DNNs $Net_{t_f}$ and $Net_{\Delta m}$ adopt MSE as the loss function which is similar to the gravity DNN in Eq.~\eqref{equ: the mean square error of gravity DNN}. Meanwhile, the DNN $Net_{\alpha}$ uses the mean absolute error (MAE) as the loss function that is defined as   
\begin{equation}
\label{equ: the mean absolute error of control DNN}
{L(\bm{\omega})} = \frac{1}{N}\sum\limits_{i = 1}^N {{{|Net_{\alpha}(\bm{x}_i|\bm{\omega}^{\alpha}) - {\alpha_i^*}|}}} 
\end{equation}
where $\alpha$ represents one of the thrust directions $\alpha_1$ , $\alpha_2$, and $\alpha_3$, and the symbol $\bm{w}^{\alpha}$ represents the parameter vector of the networks. Similar to the gravity DNN developed in Section \ref{sec: Gravity Approximation}, the algorithm Adam is employed to minimize the MSE and MAE of the five networks wherein the learning rate (LR) is set to 0.0001 and the batch number $N$ is 1000.

\begin{figure}[htbp]
\center
\vspace {2mm}
\centerline{\includegraphics[width= 6.8 in]{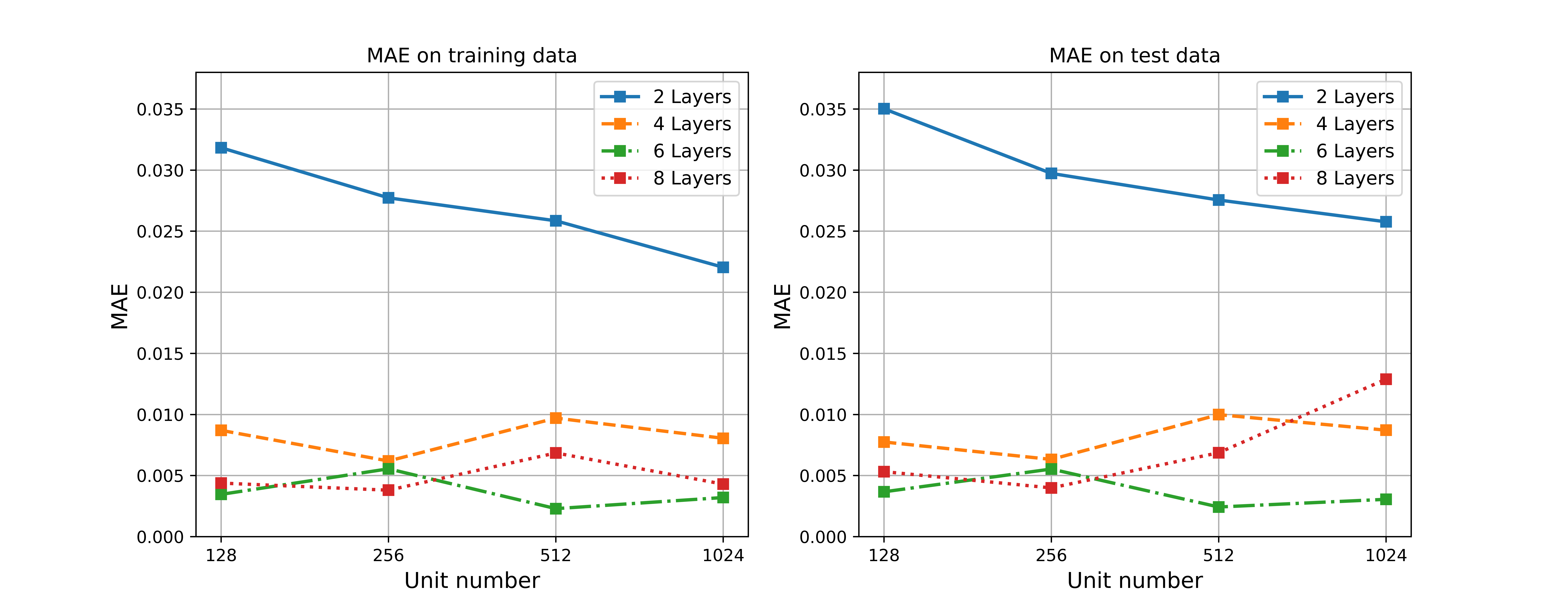}}
\vspace {-1mm}
\caption{Approximation loss comparison of Network $Net_{\alpha_1}$ with different sizes}
\label{fig: Network Size Optimization}
\vspace {-3mm}
\end{figure}

The network size is an important user-defined parameter to guarantee the quality of network approximation. In this study, the network size is optimized based on its performance on the training and test datasets. Taking the size optimization of $Net_{\alpha_1}$ as an example, Fig.~\ref{fig: Network Size Optimization} displays the MAE comparison of DNNs with different layers and units. We can obtain the following conclusions from this figure. First, compared to the increment of the network width, the network depth is more conducive to improve the learning performance of the DNNs. Second, a too mall network size would lead to the under-fitting phenomenon of the DNN, which means that the DNN can not achieve enough accuracy on the training dataset. While a too big network size also causes the deterioration of the learning effect, and these will be an obvious gap between the training error and test error, which is termed as the over-fitting phenomenon. Based on the above analysis, the network size of $Net_{\alpha}$ is determined as 6 Layers/512 Units. By the similar means, the sizes of the DNNs $Net_{t_f}$ and $Net_{\Delta m}$ are set to be 6 Layer/256 Units.

The samples of flight state versus optimal control actions are obtained by the approximate indirect method with high computational efficiency. However, there are numerical differences between the elements in both the input and output vectors of the samples, because these elements have different physical meanings and magnitudes. In this study, the Scikit-learn package of machine learning in Python environment is used to linearly convert the values of sample elements into the range of [-1,1]. This data normalization can improve the data identification ability and thus boost the prediction accuracy of the DNNs. Additionally, there is a data anti-normalization module to calculate the actual action predictions from the normalized data. Finally, the pseudo code of the training algorithm for optimal action approximation is summarized as follows in Algorithm \ref{code: Supervised training algorithm of control DNNs}, which is actually implemented in Python Environment using Tensorflow.

\begin{algorithm}[h]  
\caption{ Supervised training algorithm of the control DNNs}  
\label{code: Supervised training algorithm of control DNNs}  
\begin{algorithmic}[1]  
\STATE Randomly initialize $Net_{\alpha}(\bm{x}|\bm{w}^{\alpha})$, $ Net_{t_f} (\bm{x}|\bm{w}^{t_f})$, and $  Net_{\Delta m}(\bm{x}|\bm{w}^{\Delta m})$ with weights $\bm{w}^{\alpha}$, $\bm{w}^{t_f}$ and $\bm{w}^{\Delta m}$	
\STATE Preprogress the training data with the normalization technique \
\FOR{episode = 1, M}  
\STATE   Select a random minibatch of $N$ samples $[\bm{x}_i, \bm{\alpha}_i^*, {t_f}_i^*, \Delta m_i^*]$ from normalized training dataset 
\STATE   Update $Net_{\alpha}(\bm{x}|\bm{\omega}^{\alpha})$ with Adam algorithm by minimizing the MAE loss. 
\STATE   Update $Net_{t_f}(\bm{x}|\bm{\omega}^{t_f})$ and $Net_{\Delta m}(\bm{x}|\bm{\omega}^{\Delta m})$  with Adam algorithm by minimizing the MSE loss. 
\ENDFOR 
\end{algorithmic}  
\end{algorithm}

\subsection{DNN-Based Real-Time Optimal Landing}

\begin{figure}[htbp]
\center
\vspace {2mm}
\centerline{\includegraphics[width= 5 in]{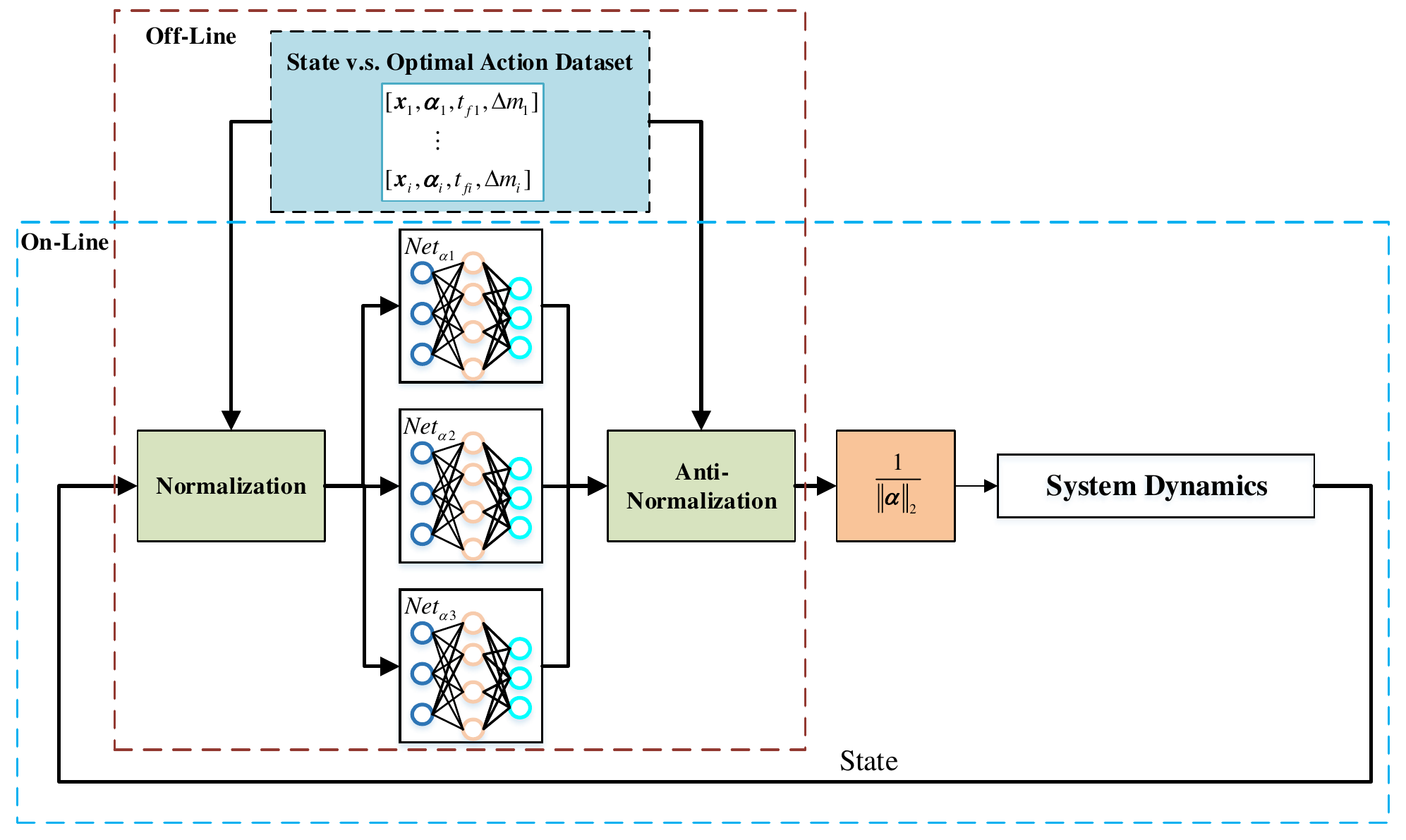}}
\vspace {-1mm}
\caption{DNN-based real-time optimal controller}
\label{fig: DNN-based real-time optimal controller}
\vspace {-3mm}
\end{figure}

Through learning the training data generated by the proposed approximate indirect method, the DNNs are trained to approximate the state-based optimal actions. Figure \ref{fig: DNN-based real-time optimal controller} describes the real-time optimal control for asteroid landing driven by the trained DNNs. Considering the preprocessing of the training data, a normalization module is used for network inputs and an anti-normalization module is used for network predictions. Additionally, a multiplier of $1/\Vert \bm{\alpha} \Vert_2 $ is introduced to normalize the thrust direction vector $\bm{\alpha}$. Different from the traditional indirect and direct methods, the well-trained DNNs can directly predict the optimal control actions according to the flight state and no longer needs to solve the original OCPs onboard. Consequently, the DNN-based landing controller is capable of real-time optimal control for asteroid landing missions. Multiple experiments will be carried out to evaluate the performance of this DNN-based landing controller in the next section.

\section{Simulation Results}
\label{sec: Simulation and Results}

Numerical simulations are presented in this section to evaluate the effectiveness and performance of the proposed DNN-based landing control for asteroids with irregular gravitational fields. The evaluation of the DNN-based landing controller involves several aspects including the real-time performance of command generation, the accuracy of optimal action prediction, and the generalization capability with respect to the states. Furthermore, to quantitatively evaluate the proposed control process, the solution optimality, terminal landing accuracy, and control robustness under uncertainties are also investigated. Based on these considerations, four simulation experiments are conducted, and the performance analysis of the controller is discussed in detail. The proposed methodology has been implemented in C and Python on a desktop with Inter Core i7-7700 CPU @3.60GHz RAM 8.00GB and operating system of Windows 10. The near-Earth asteroid 443 Eros is used in the considered landing scenarios.


\subsection{Real-Time Performance}






One of the biggest advantages of the proposed method is its real-time performance for landing control using DNNs. Specifically, the gravity DNN developed in Section \ref{sec: Gravity Approximation} is used to improve the speed of trajectory optimization and the control DNNs developed in Section \ref{sec: DNN-Based Real-Time Optimal for Landing} are developed to achieve real-time decision-makings for on-board landing control. Considering that the polyhedral method is implemented in C environment while our DNNs are developed in Tensorflow using Python environment, we have to evaluate the computational efficiency of these two simulation environments first. 

Table \ref{tab: Time Comparison} shows that the computational efficiency ratio between C and Python is 36:1, which is obtained through time-consumption analysis of a series of benchmark algorithms. The computational efficiency comparison of these two environments is complicated and may be influenced by the running algorithms, and the software and hardware of the computers. However, the running of algorithms in C is generally faster than in Python. The awareness of the above conclusion does not affect the results of the subsequent time comparisons, however, it will be beneficial for further studying the efficiency of the DNN-based controller.

\begin{table}[htb]
\centering  
\caption{ Time cost comparison}
\label{tab: Time Comparison}
\begin{tabular}{ccccc}
 \hline
 \hline
                                         & \multicolumn{2}{c}{Polyhedron (C)} & \multicolumn{2}{c}{DNN (Python)} \\
Computation Efficiency Ratio             & \multicolumn{2}{c}{36}             & \multicolumn{2}{c}{1} \\
Gravity Calculation                      & \multicolumn{2}{c}{41.3 ms}        & \multicolumn{2}{c}{0.8 ms}\\
Trajectory Optimization without Homotopy & \multicolumn{2}{c}{295.2 min}      & \multicolumn{2}{c}{5.7 min}\\
Trajectory Optimization with Homotopy    & \multicolumn{2}{c}{48.5 min}       & \multicolumn{2}{c}{1.1 min} \\
 \hline
 \hline
\end{tabular}
\end{table}

First, the speed improvement in trajectory optimization benefited from the gravity DNN is investigated. A polyhedral model with 25,350 vertices and 49,152 faces is used for 433 Eros, and its calculation is time-consuming. As displayed in Table \ref{tab: Time Comparison}, the gravity calculation using the polyhedral method consumes about 41.3 ms in C. Meanwhile, the well-trained gravity DNN $Net_{\mathrm {G}}(\bm{r}|\bm{\omega})$ can predict the accurate value of gravity in 0.8 ms. As such, leveraging the DNN technique, the calculation speed of gravity is improved by 51 times (without considering the impact of simulation environment), and the time consumption of trajectory optimization is reduced from 295.2 min to 5.7 min. Moreover, the homotopy technique can further decrease the time consumption by improving the success rate of the approximate indirect method. Considering that the running time of algorithms in C is usually far less than that in Python, we anticipate that the implementing the DNN in C would bring much more speed increment to trajectory optimization. 

In addition, the real-time performance is also observed in control command generation by the developed DNN-based landing controller, which is capable of directly predicting the optimal control instructions based on the flight state. Without the need to solve the OCPs onboard, the optimal control commands can be generated within 0.9 ms with a much higher updating frequency than any direct or indirect method. Additionally, the DNN-based trajectory controller does not suffer the convergence issues and can provide deterministic solutions to asteroid landing problems. The real-time performance and stable calculation process of the DNN-based landing controller show great potential for practical engineering applications \cite{sanchez2018real}.

\subsection{Optimality Analysis}

\begin{figure}[htbp]
\center
\vspace {2mm}
\centerline{\includegraphics[width= 7 in]{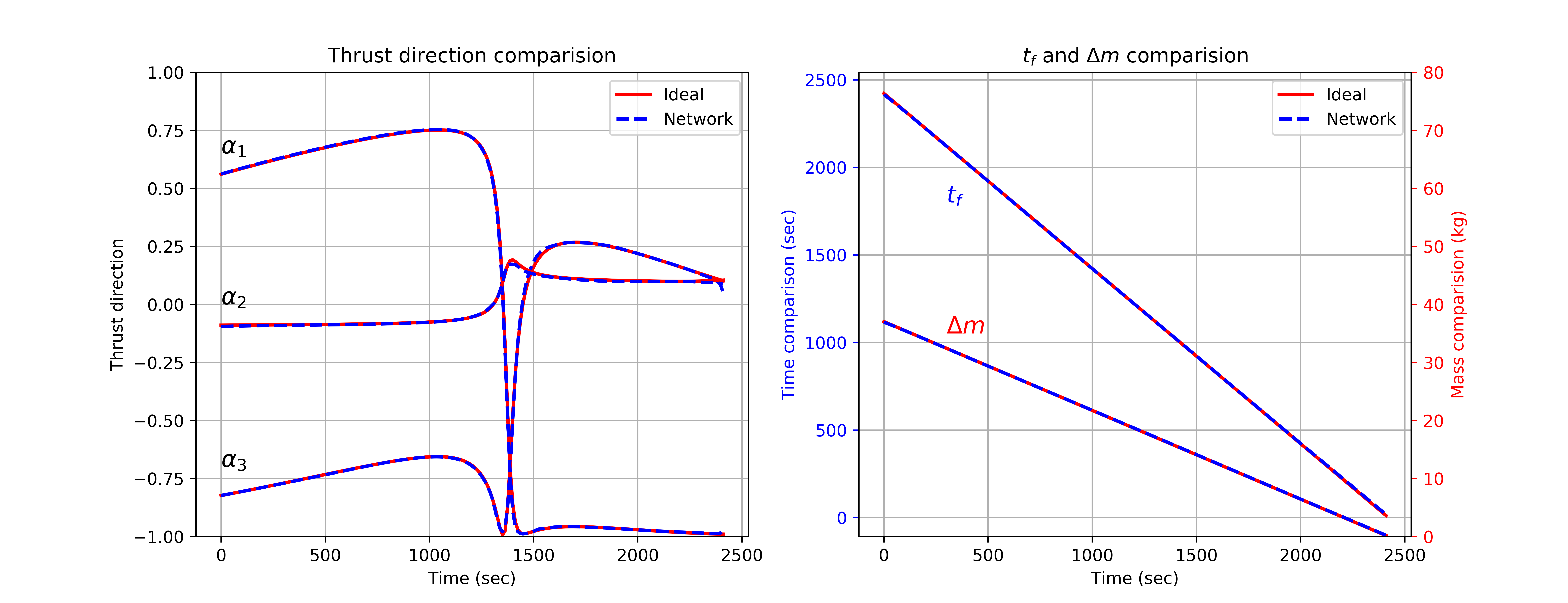}}
\vspace {-1mm}
\caption{Comparison of optimal actions between ideal solution and DNN prediction}
\label{fig: DNN prediction comparison on optimal actions}
\vspace {-3mm}
\end{figure}

The network approximation accuracy to the optimal actions and the state-based generalization ability are two factors that determine the performance of the developed DNN-based landing controller and are needed to be investigated. Figure \ref{fig: DNN prediction comparison on optimal actions} displays the approximation results of the DNNs in learning the thrust directions $\bm{\alpha}$, time-to-go $t_f$, and mass consumption $\Delta m$. In this figure, all red solid lines represent the ideal outputs referring to the training dataset, while the blue dotted lines indicate the predictions of the trained DNNs. The left subgraph displays the approximation results to the three thrust direction components and the right one depicts the approximation results to the time-to-go $t_f$ and the mass consumption $\Delta m$. The right subgraph has two Y-axes, one of which indicates the time-to-go and the other one denotes the mass consumption. As we can see from the two subgraphs, the approximate optimal solutions obtained by the trained DNNs are very accurate. Even at the state with sharp changes in the instructions, the DNN results can still follow the ideal solutions very well. Additionally, the functional relationships of the time and mass consumption with respect to the states are relatively simple, thus, their approximations seem much better than the approximation of the thrust directions.

\begin{figure}[htbp]
\center
\vspace {2mm}
\centerline{\includegraphics[width= 7 in]{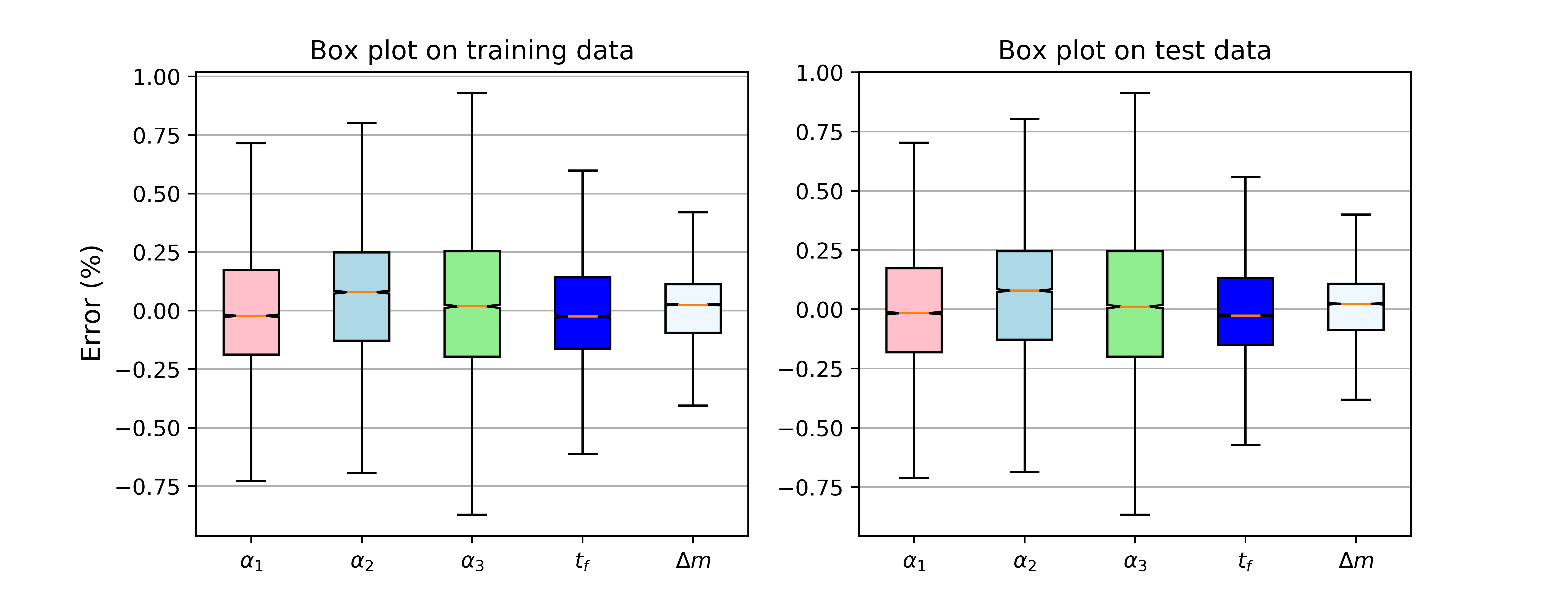}}
\vspace {-1mm}
\caption{Box plot of DNN prediction errors on optimal actions}
\label{fig: Box plot of control network}
\vspace {-3mm}
\end{figure}

Then, we further quantitatively analyze the network approximation errors using the box plot. Figure \ref{fig: Box plot of control network} shows the statistical results of approximation errors of the five control DNNs. This box plot is drawn using 10,000 samples randomly chosen from the training dataset and the test dataset. The Errors are calculated according to Eq.~\eqref{equ: relative error vector}. As we can see from the figure that the relative approximation errors of the three thrust directions stay within $\pm 0.25\%$ and the maxima of their absolute values are smaller than $1\%$. Due to the relatively simple functional relationships of the time and mass consumption with respect to the states, the approximation errors of the DNNs $Net_{t_f}$ and $Net_{\Delta m}$ are much smaller than other networks with the means within $ \pm 0.15\%$ and the limitations within $\pm 0.6\%$. Based on these analysis, we can draw the conclusion that the trained control DNNs have high approximate accuracy to the optimal actions and are potential for on-board solution generations for asteroid landing missions.

\begin{figure}[htbp]
\center
\begin{minipage}[t]{0.48\linewidth}
\centerline{\includegraphics[width= 3.1 in]{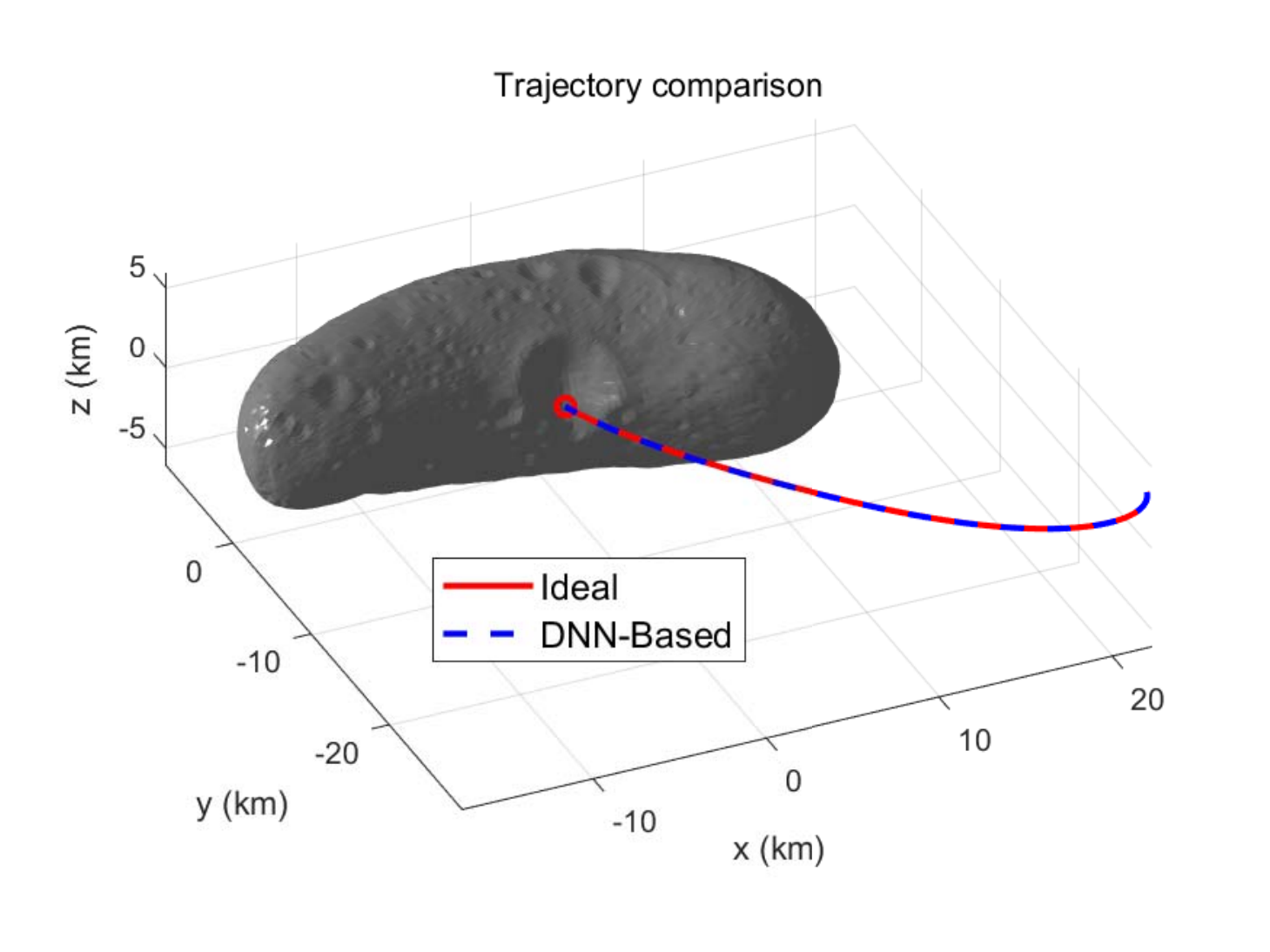}}
\centerline{\textbf{(a)}}
\end{minipage}
\begin{minipage}[t]{0.48\linewidth}
\centerline{\includegraphics[width= 3 in]{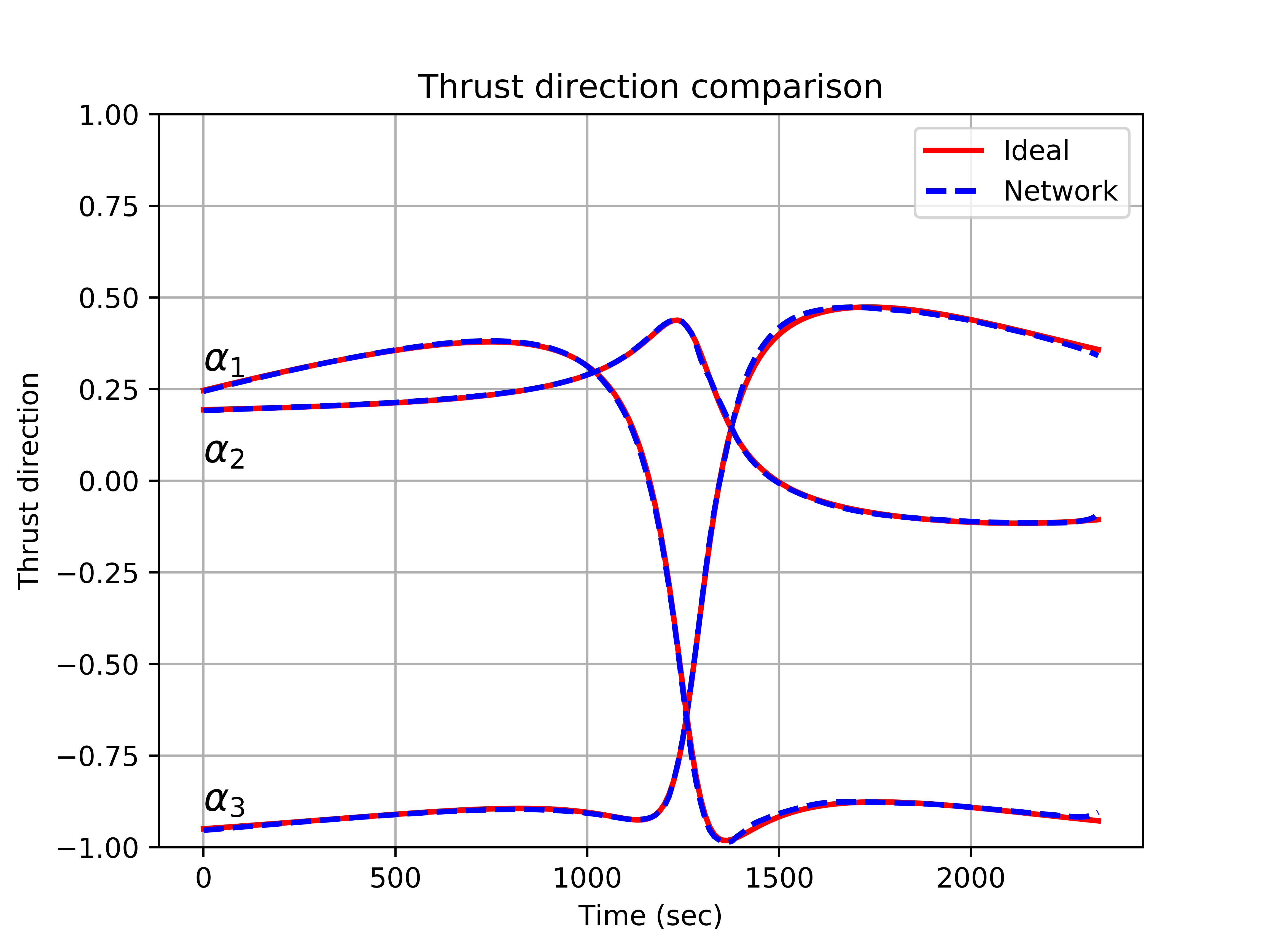}}
\centerline{\textbf{(b)}}
\end{minipage}
\caption{Case1: Solutions driven by the DNN-based landing controller  where (a) shows the trajectories and (b) shows the controls}
\label{fig: Landing driven by the DNN-based trajectory controller (case1)}
\end{figure}

\begin{figure}[htbp]
\center
\begin{minipage}[t]{0.48\linewidth}
\centerline{\includegraphics[width= 3.1 in]{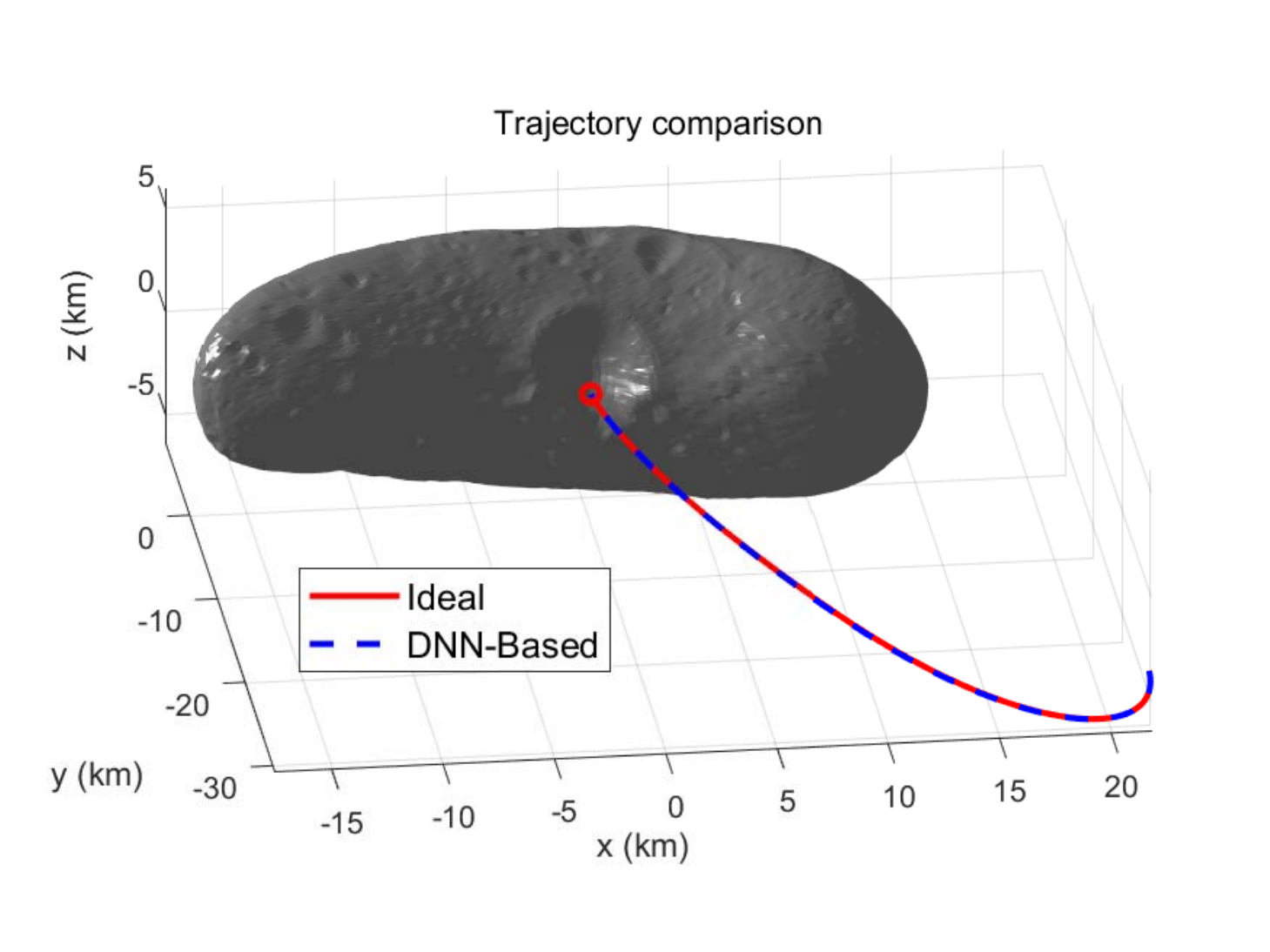}}
\centerline{\textbf{(a)}}
\end{minipage}
\begin{minipage}[t]{0.48\linewidth}
\centerline{\includegraphics[width= 3 in]{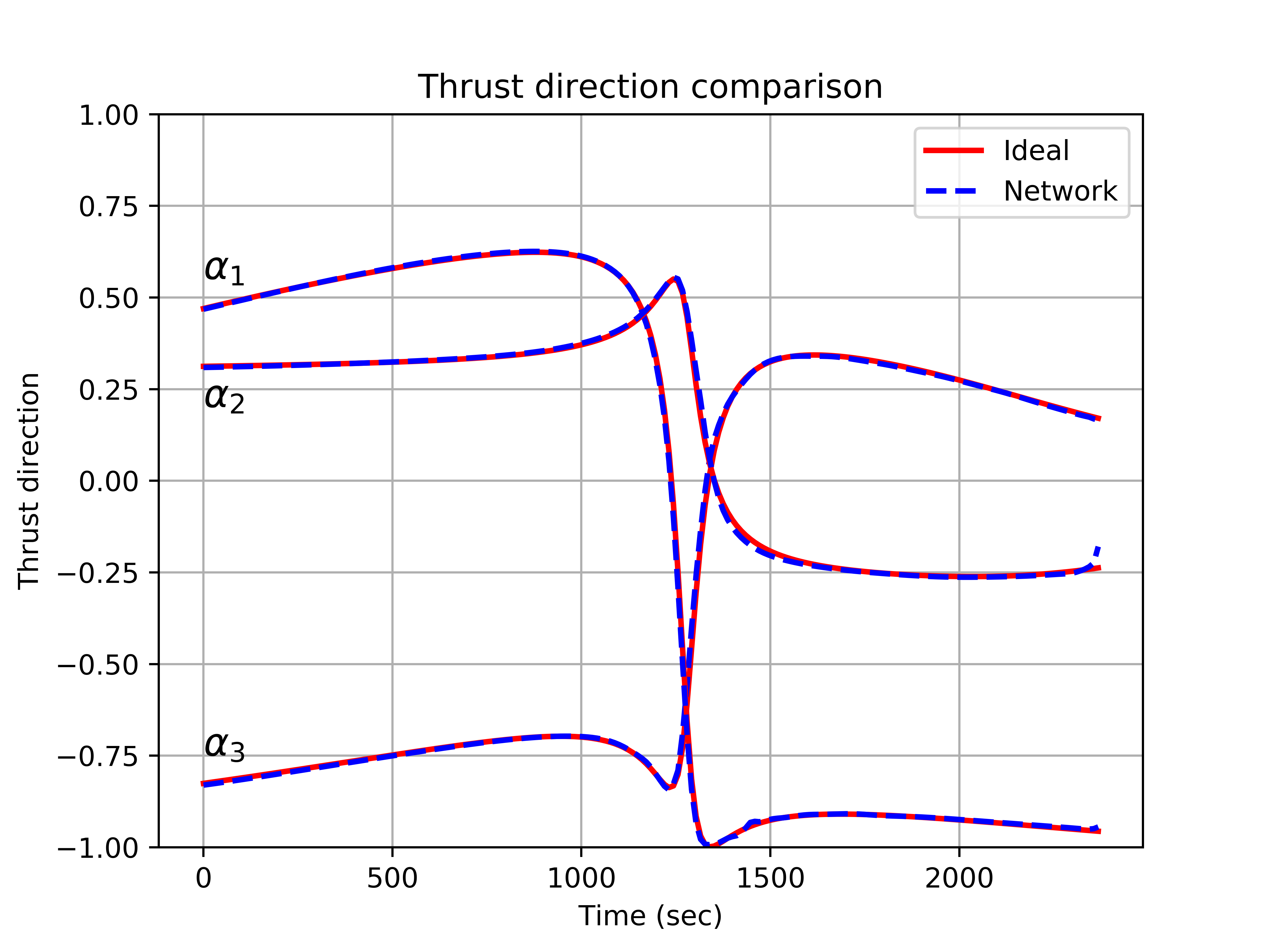}}
\centerline{\textbf{(b)}}
\end{minipage}
\caption{Case2: Solutions driven by the DNN-based landing controller where (a) shows the trajectories and (b) shows the controls}
\label{fig: Landing driven by the DNN-based trajectory controller (case2)}
\end{figure}

In addition, the asteroid landing flight is a long-lasting control process, and small network fitting error may propagate and cause large flight error. However, large trajectory deviation can be reduced through continuous command corrections. Considering the existing approximation errors of DNNs, two flight cases driven by the DNN-based landing controller are carried out to test the controller's performance on reducing the propagated errors in a closed-loop control manner. The lander's initial state is randomly selected from the initial space illustrated in Fig.~\ref{fig: DataGeneration}. The whole flight is completely driven by the DNN-based landing controller, which independently determines the control instructions according to the actual flight state. The simulation is terminated based on the predicted landing time from the network $Net_{t_f}$. The resulting flight profiles are compared to the ideal flight solutions obtained by the indirect method and shown in Figs.~\ref{fig: Landing driven by the DNN-based trajectory controller (case1)} and \ref{fig: Landing driven by the DNN-based trajectory controller (case2)}. In these two figures, the left subgraphs exhibit the trajectories from these two approaches and the right subgraphs compare the corresponding profiles of the three thrust directions $\alpha_1$, $\alpha_2$, and $\alpha_3$. As it can be seen from the figures that the developed DNN-based landing controller is capable of independently driving the lander to the target state in an optimal manner, and the driven trajectories coincide with the ideal trajectories designed off-line by the indirect method. The above two experiments further verify the optimality of the solution generated by the proposed DNN-based controller for asteroid landing missions. 



\subsection{Terminal Guidance Error}

In addition to the real-time calculation process and the guaranteed solution optimality, the accuracy at the terminal landing moment is also critical to the success of the entire mission. In this subsection, we will discuss the terminal guidance accuracy of the asteroid landing flight  driven by the proposed DNN-based landing controller. 


\begin{figure}[htbp]
\center
\vspace {2mm}
\centerline{\includegraphics[width= 4.5 in]{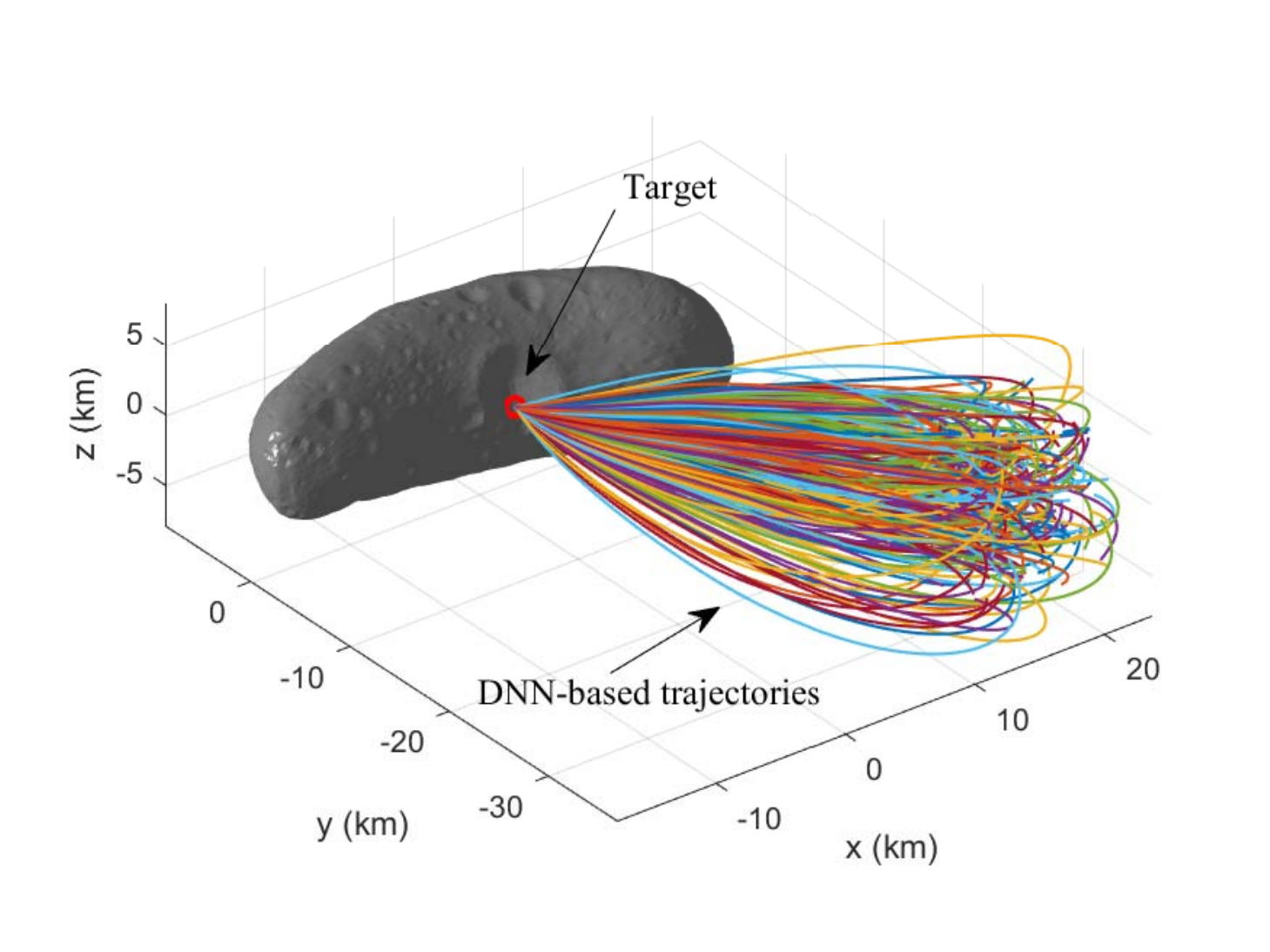}}
\vspace {-1mm}
\caption{DNN-based landing trajectories}
\label{fig: DNN-based landing trajectories}
\vspace {-3mm}
\end{figure}

\begin{figure}[htbp]
\center
\begin{minipage}[t]{0.45\linewidth}
\centerline{\includegraphics[width= 3.2 in]{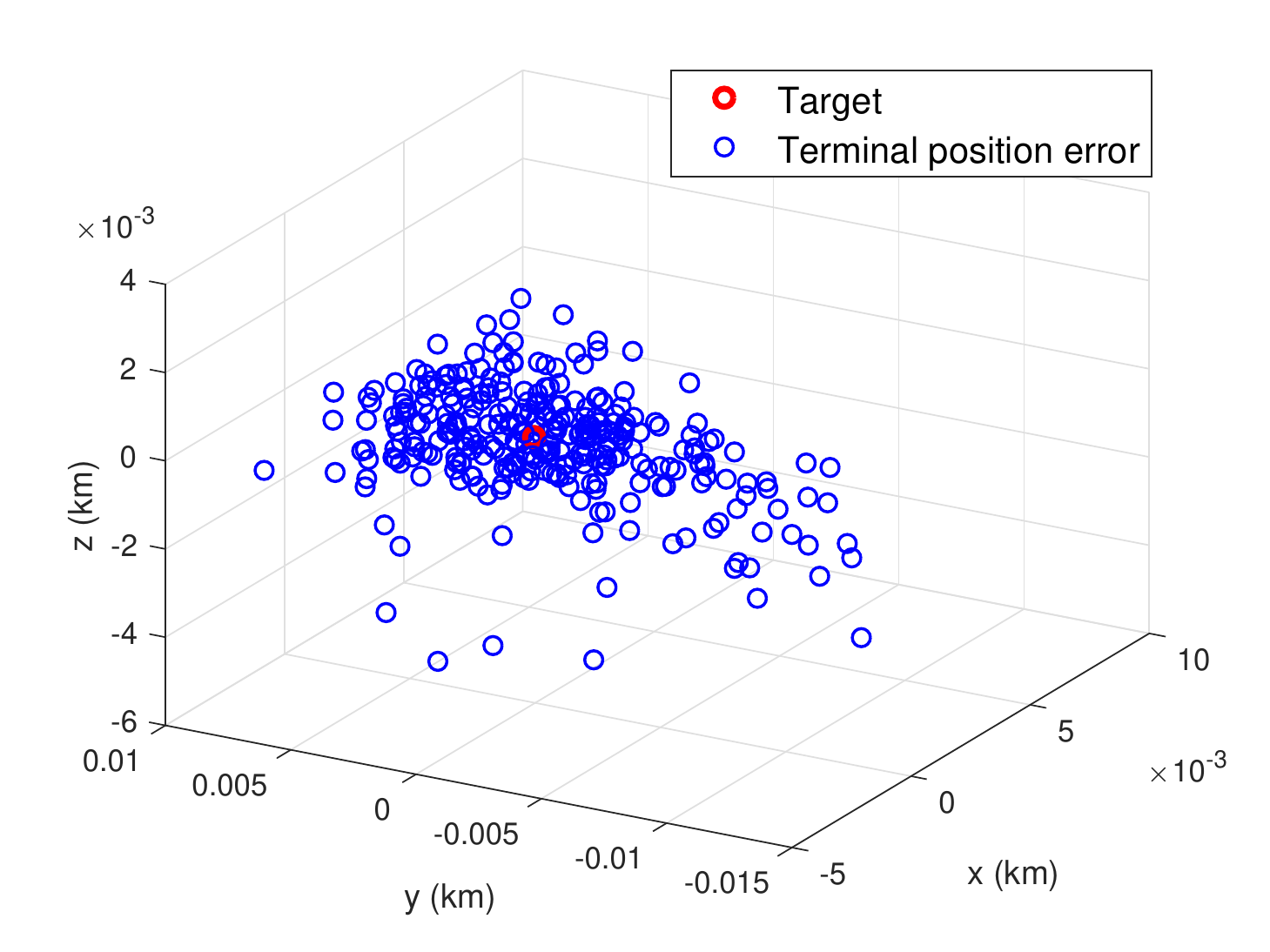}}
\centerline{\textbf{(a)}}
\end{minipage}
\begin{minipage}[t]{0.5\linewidth}
\centerline{\includegraphics[width= 3.2 in]{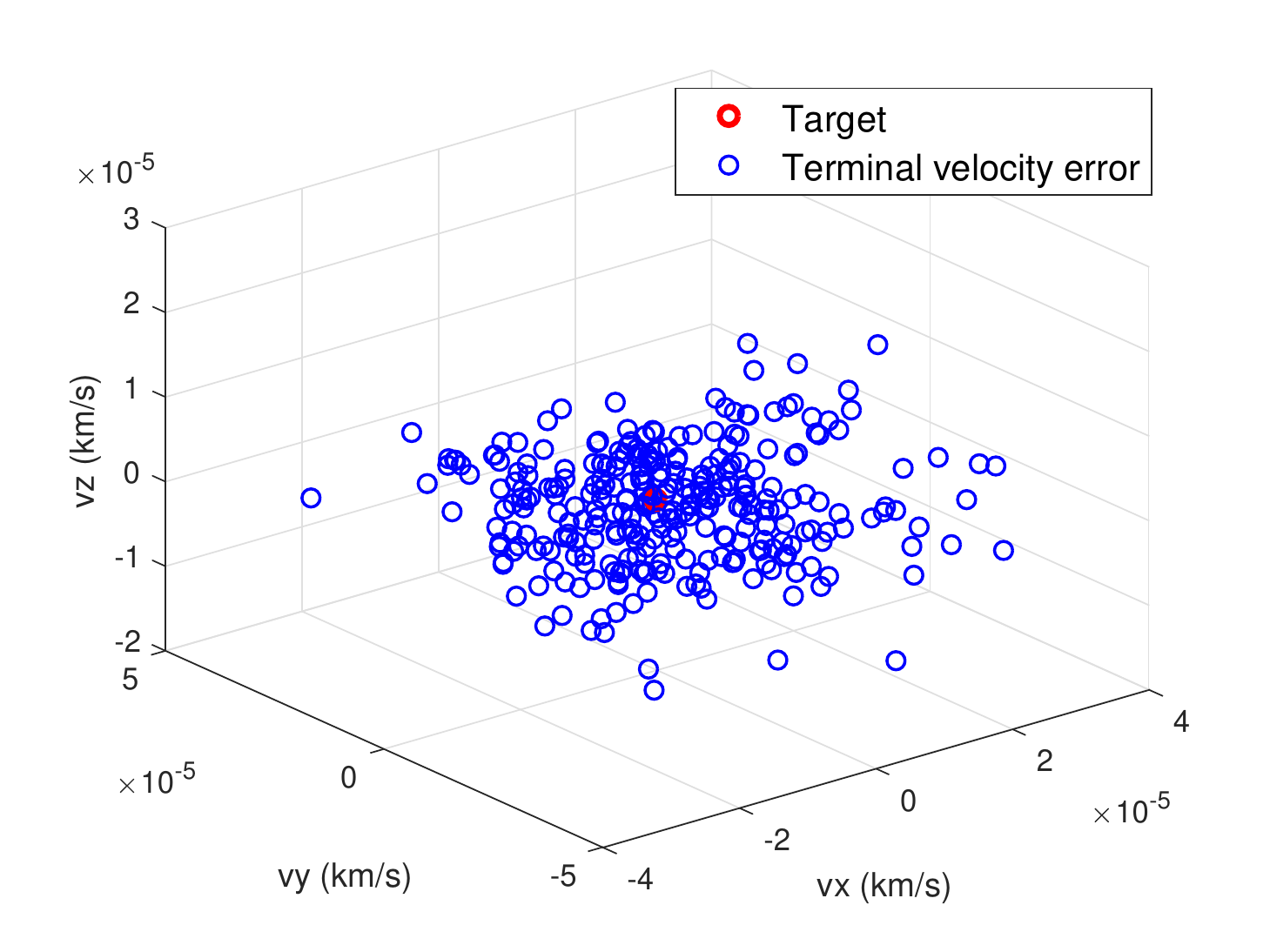}}
\centerline{\textbf{(b)}}
\end{minipage}
\caption{Terminal guidance errors driven by the DNN-based landing controller where (a) shows the terminal position errors and (b) shows the terminal velocity errors}
\label{fig: Terminal guidance error driven by the DNN-based landing controller}
\end{figure}

Figure \ref{fig: DNN-based landing trajectories} displays 300 DNN-driven trajectories starting from random initial states. As it is shown in the figure that although the start conditions of the lander are randomly initialized, the DNN-based landing controller can autonomously decide the optimal control instructions according to the flight states and steer the lander to the target site accurately. Considering that the flight state is continuously changing and there must be states that do not exist in the training dataset, we can conclude that this DNN-based controller has excellent generalization ability for state-based decision-makings. Figure \ref{fig: Terminal guidance error driven by the DNN-based landing controller} illustrates the terminal positions and velocities of these 300 trajectories. We can see from the figure that the terminal guidance position errors are within the range of $\pm 0.01$km, and the terminal velocity errors are within the range of $\pm 4\times 10^{-5}$ km/s. Hence, we can conclude that the DNN-based landing controller is capable of steering the lander to the target point with high accuracy. The terminal guidance accuracy can be further improved using some auxiliary techniques, for example, a multi-scale network cooperation strategy was proposed in \cite{chenglin2018TAES} to improve the terminal flight accuracy of solar sail orbit transfers. In addition, a compound control scheme was developed in \cite{chenglin2018RealTimeAstrodynamics} to improve the terminal landing accuracy for the Moon landing mission. 


\subsection{Robustness Analysis}

\begin{figure}[htbp]
\center
\vspace {2mm}
\centerline{\includegraphics[width= 7 in]{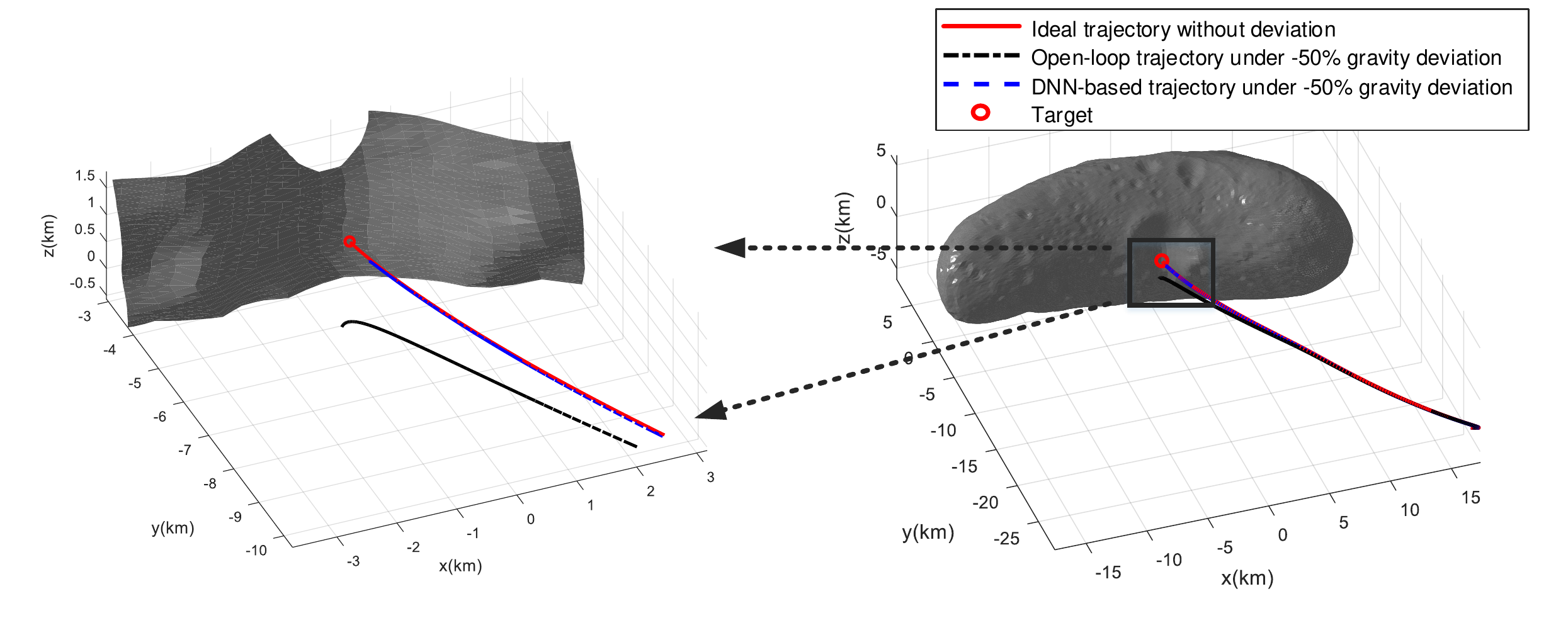}}
\vspace {-1mm}
\caption{DNN-based driving trajectory with a $-50\%$ gravity deviation}
\label{fig: Robuseness}
\vspace {-3mm}
\end{figure}

In consideration of the inevitable dynamical uncertainties caused by the approximation error of gravitational force, the solar wind, and the aberration of light, we test the robustness of the DNN-based landing controller in this subsection. For preliminary analysis, we assume that the gravity calculation errors and all other force interferences are considered together as the deviation of the gravitational field. Figure \ref{fig: Robuseness} shown the trajectories under the gravity with $-50\%$ deviation. The red solid trajectory is the ideal solution without gravity deviations. Under the effect of gravity deviation, the black dotted trajectory is obtained by the open-loop indirect method and the blue dashed trajectory is the solution by the DNN-based controller. As it can be seen from the figure that the flight by the open-loop indirect method has big terminal state error. Traditionally, this error could be reduced by a closed-loop trajectory tracking system. Designing a closed-loop tracking scheme is not an easy task due to the difficulty of generating a reference trajectory in advance, the risk of trajectory tracking failure, and the lack of flight autonomy. In this study, the proposed DNN-based landing controller can adjust the control instruction on-line based on actual flight state in every guidance cycle. Consequently, the proposed DNN-based controller can significantly reduce the flight errors through rolling command corrections and achieve high model robustness and terminal landing accuracy, even when large deviations occur in the gravity field.

\begin{table}[htb]
\centering  
\caption{ Robustness analysis}
\label{tab: Robustness}
\begin{tabular}{c|cccc|cccc}
 \hline
 \hline
 & \multicolumn{4}{c|}{Indirect Method} & \multicolumn{4}{c}{DNN-Based Control} \\
  \cline{2-9}
 & \multicolumn{2}{c}{Position Error (km)} & \multicolumn{2}{c|}{Velocity Error (km/s)} & \multicolumn{2}{c}{Position Error (km)} & \multicolumn{2}{c}{Velocity Error (km/s)} \\
  \cline{2-9}
          & Mean  & Max   & Mean     & \multicolumn{1}{c|}{Max}      & Mean    & Max    & Mean     & Max  \\
 \hline
 $-100\%$ & 3.37    & 5.28    & 4.98E-3 & 6.43E-3 & 1.25E-1 & 1.37E-1 & 2.33E-4  & 3.96E-4  \\
 $-50\%$  & 1.75    & 2.78    & 2.66E-3 & 3.52E-3 & 7.72E-2 & 8.62E-2 & 6.24E-5 & 1.88E-4  \\
 $-10\%$  & 3.61E-1 & 5.81E-1 & 5.59E-4 & 7.51E-4 & 1.81E-2 & 2.73E-2 & 1.08E-4  & 1.37E-4  \\
 $0$      & 0       & 0       & 0       & 0       & 3.29E-3 & 1.10E-2 & 1.51E-5  & 3.93E-5  \\
 $+10\%$  & 3.66E-1 & 5.91E-1 & 5.53E-4 & 7.38E-4 & 1.52E-2 & 1.97E-2 & 9.71E-5  & 1.21E-4    \\
 $+50\%$  & 1.735   & 2.80    & 1.66E-3 & 2.45E-3 & 8.36E-2 & 9.11E-2 & 5.08E-4  & 5.58E-4    \\
 $+100\%$ & 3.405   & 5.71    & 3.28E-3 & 5.20E-3 & 1.61E-1 & 1.77E-1 & 6.32E-4  & 7.02E-4    \\
 \hline
 \hline
\end{tabular}
\end{table}

Table \ref{tab: Robustness} presents the statistical results of terminal flight errors under  different gravity deviations. As we can see from the table that the open-loop control scheme based on the indirect method has relatively large flight errors and behaves poor robustness when obvious uncertainties are incorporated in the dynamics. In contrast, benefiting from its rolling trajectory correction capability, the DNN-based landing controller is more robust and achieves more accurate terminal flight accuracy than the traditional indirect methods. 

\section{Conclusion}
\label{sec: conclusion}

In this study, deep neural network technology is introduced to improve the real-time performance and flight autonomy for asteroid landing control. This work processes the following three contributions. First, a new DNN-based gravity model is developed to replace the time-consuming homogeneous polyhedral method and significantly improve the computational efficiency of gravity determination. Then, an improved approximate indirect method is proposed to solve the time-optimal landing problems where the gravity term in the dynamical system is replaced with the designed DNN-based gravity model and a homotopy method is employed to further reduce the computational cost by connecting the original time-optimal problem to a gravity-free time-optimal problem. Furthermore, five DNNs are developed to learn the state-based optimal actions obtained from the approximate indirect method, and a DNN-based landing controller composed of these DNNs is eventually designed to successfully achieve the real-time optimal control for asteroid landing missions. Simulation results of time-optimal landing problems of landing on 443 Eros are given to substantiate the improvement of real-time performance and illustrate the effectiveness of the developed DNN-based controller in solution optimality, terminal flight accuracy, and robustness under uncertainties.

\section*{Acknowledgments}
This work is supported by the National Natural Science Foundation of China (Grants 11872223, 11672146 and 11432001).

\bibliographystyle{unsrt}
\bibliography{myreference} 

\end{document}